\documentclass[a4paper,11pt]{article}

\usepackage[small,bf,hang]{caption} 

\usepackage{todonotes}
\usepackage{amsmath,amsthm,amssymb}
\usepackage{graphicx,subcaption,tikz}
\usepackage{comment,xspace,paralist}							
\usepackage[shortlabels]{enumitem}	
\usepackage{hyperref}
\usepackage{changepage}
\hypersetup{
    colorlinks=true,
    citecolor = blue,
    linkcolor=blue,
    filecolor=magenta,
    urlcolor=cyan,
}
\usepackage[margin=1in]{geometry}
\usepackage{verbatim}
\usepackage[capitalize]{cleveref}
\usepackage{color}
\usepackage{authblk}

\usetikzlibrary{decorations.pathreplacing, calligraphy}

\newtheorem{theorem}{Theorem}[section]
\newtheorem{lemma}[theorem]{Lemma}

\newtheorem{claim}{Claim}[theorem]
\Crefname{claim}{Claim}{Claims}

\Crefname{subsection}{Subsection}{Subsections}
\newtheorem{conjecture}{Conjecture}
\Crefname{conjecture}{Conjecture}{Conjectures}
\theoremstyle{definition}
\newtheorem{definition}[theorem]{Definition}

\expandafter\let\expandafter\oldproof\csname\string\proof\endcsname
\let\oldendproof\endproof
\renewenvironment{proof}[1][\proofname]{%
	\oldproof[\normalfont\bfseries #1]%
}{\oldendproof}
\newenvironment{subproof}[1][\normalfont\it Subproof]{%
	\begin{proof}[#1]%
	}{%
	\end{proof}%
}

\newcommand{\eps}{\varepsilon}
\newcommand{\EH}{\mathrm{EH}}
\newcommand{\cG}{\mathcal G}
\newcommand{\cC}{\mathcal C}
\newcommand{\cO}{\mathcal O}
\newcommand{\Rs}{\mathsf R}
\newcommand{\G}{\overline{G}}
\DeclareMathOperator{\homog}{hom}

\title{The Erd\H{o}s--Hajnal property for paths}

	\author[1]{Shenwei Huang}

    \author[2]{Yidong Zhou}

    \author[3]{Yuwen Zhou}
	
	\affil[1] {School of Mathematical Sciences and LPMC, Nankai University, Tianjin 300071, P.R.~China. Email: {\tt shenweihuang@nankai.edu.cn}.}
  
    \affil[2] {College of Computer Science, Nankai University, Tianjin 300350, China.}

    \affil[3] {Irving K. Barber Faculty of Science, University of British Columbia Okanagan, Kelowna, BC V1V 1V7, Canada.}

\begin{document}

\maketitle

\begin{abstract}
    The well-known Erd\H{o}s-Hajnal conjecture states that for any graph $H$, there is a constant $c=c(H)>0$ such that every $n$-vertex graph $G$ with no induced copies of $H$ contains a clique or an independent set of size at least $n^{c}$.
    In this paper, we prove that the Erd\H{o}s-Hajnal conjecture holds for every path. 
\end{abstract}

\section{Introduction}\label{sec:intro}

    A {\em clique} ({\em stable set}) is a vertex set whose vertices are pairwise adjacent (nonadjacent, resp.).
    Ramsey~\cite{Ra30} proved that for every integer $k\geq 1$, there exists a positive integer $R(k,k)$ such that every graph on at least $R(k,k)$ vertices contains a clique or stable set of size $k$. 
    This shows that a sufficiently large graph must contain a large clique or stable set. 
    Erd\H{o}s \cite{Er47} proved via the probabilistic method that for large $n$, there exists an $n$-vertex graph with no cliques or stable sets of size at least $2\log n$. 
    This shows that general graphs do not necessarily contain a clique or stable set of polynomial size. 
    However, $H$-free graphs might behave differently.

    Let $G,H$ be graphs. 
    We say that $G$ is {\em $H$-free} if $G$ has no induced subgraph that is isomorphic to $H$. 
    In 1989, Erd\H{o}s and Hajnal \cite{eh} showed that for any graph $H$, there exists a constant $c=c(H)>0$ such that any $n$-vertex $H$-free graph has a clique or stable set of size at least $e^{c\sqrt{\log n}}$.
    Erd\H{o}s and Hajnal~\cite{eh} further conjectured that one can improve the bound to a polynomial of $n$, which has become known as the Erd\H{o}s-Hajnal conjecture.

    \begin{conjecture}\label{conj:E-H}
        \textnormal{(Erd\H{o}s-Hajnal \cite{eh})} For every graph $H$, there exists $c=c(H)>0$ such that every $n$-vertex $H$-free graph has a clique or stable set of size at least $n^c$.
    \end{conjecture}

    We say a graph $H$ has the {\em Erd\H{o}s-Hajnal property} if $H$ satisfies Conjecture \ref{conj:E-H}. 
    The property is preserved under complementation and, by a theorem of Alon, Pach and Solymosi \cite{aps}, under substitution. 
    It holds for all graphs with at most four vertices, for the bull \cite{cs}, for the five-cycle \cite{csss} and for the five-vertex path $P_5$ \cite{nss7}. 
    The last result completed the case of all five-vertex graphs. Huang, Ju and Zhou \cite{hjz} developed the method of \cite{nss7} into a general framework and used it for two six-vertex graphs. 
    We refer to \cite{chudnovsky} for a survey.

    Paths are a natural family to consider next. 
    Nguyen, Scott and Seymour \cite{nss5} proved that every $P_s$-free graph on $n$ vertices has a clique or stable set of size $2^{(\log n)^{1-o(1)}}$. 
    We prove the conjecture for all paths.

    \begin{theorem}\label{thm:main}
        For every integer $s\ge1$ there is $\tau_s>0$ such that every $P_s$-free graph $G$ has a clique or a stable set of size at least $|G|^{\tau_s}$.
    \end{theorem}

    \paragraph{Independent work.} The proof of Theorem~\ref{thm:main}, together with its Lean formalization, was made public on 7 October 2026 and archived on Zenodo on 9 October 2026 \cite{zenodo}. Later, on 17 October 2026, Nguyen \cite{nguyen26} posted on arXiv an independent proof of the Erd\H{o}s--Hajnal property for all paths, by a different method. As reported in \cite{nguyen26}, a team at Google Research and Google DeepMind has also announced an independent proof.

\section{Preliminaries}\label{sec:prelim}

    \subsection{Notation}
        All graphs are finite and simple.  
        For a graph $G$ we write $|G|$ for its number of vertices, $\G$ for its complement and $G[S]$ for the subgraph induced on $S\subseteq V(G)$. 
        We write $N(v)$ for the set of neighbors of $v$ in $G$. $P_r$ is the path with $r$ vertices. 
        A class of graphs is \emph{hereditary} if it is closed under taking induced subgraphs. 
        A hereditary class $\cC$ has the \emph{Erd\H{o}s--Hajnal property} if there is $\tau>0$ such that every $G\in\cC$ has a clique or a stable set of size at least $|G|^\tau$. 
        We call $\tau$ an \emph{Erd\H{o}s--Hajnal exponent} of $\cC$, and for a graph $H$ we write $\EH(H)$ for the statement that the class of $H$-free graphs has this property. 
        We write $\homog(G)$ for the largest size of a clique or a stable set of $G$ and $[k]=\{1,\ldots,k\}$ for an integer $k\ge1$.

        Let $A,B\subseteq V(G)$ be disjoint. 
        $A$ is \emph{complete} (\emph{anticomplete}) to $B$ if every vertex of $A$ is adjacent (nonadjacent) to every vertex of $B$. 
        A vertex $v\notin A$ is \emph{mixed} on $A$ if it has a neighbor and a nonneighbor in $A$. 
        For $x>0$, \emph{$B$ is $x$-sparse to $A$} if every vertex of $B$ has at most $x|A|$ neighbors in $A$, and $(A,B)$ is \emph{weakly $x$-sparse} if the number of edges between $A$ and $B$ is at most $x|A||B|$. 
        If $B$ is $x$-sparse to $A$, then $(A,B)$ is weakly $x$-sparse.

        $G$ is \emph{$x$-sparse} if its maximum degree is at most $x|G|$, and \emph{$x$-restricted} if $G$ or $\G$ is $x$-sparse. 
        A set $S\subseteq V(G)$ is $x$-sparse ($x$-restricted) if $G[S]$ is. 
        A graph is \emph{anticonnected} if its complement is connected, and the \emph{anticonnected components} of $G$ are the subgraphs induced on the vertex sets of the components of $\G$.

        A \emph{blockade} in $G$ is a sequence $\mathcal B=(B_1,\ldots,B_k)$ of disjoint subsets of $V(G)$.
        Its \emph{length} is $k$ and its \emph{width} is $\min_i|B_i|$. 
        It is a \emph{$(k,w)$-blockade} if its length is at least $k$ and its width is at least $w$. 
        It is \emph{pure} if every pair $(B_i,B_j)$ is complete or anticomplete, and \emph{complete} (\emph{anticomplete}) if every pair is complete (anticomplete). 
        We write $V(\mathcal B)=\bigcup_{i=1}^k B_i$.
        We call a blockade $\mathcal B$
        \begin{itemize}
            \item \emph{$x$-semisparse} if for all $i<j$, either $B_j$ is complete to $B_i$, or $B_j$ is $x$-sparse to $B_i$;
            \item \emph{weakly $x$-semisparse} if for all distinct $i,j$ the pair $(B_i,B_j)$ is complete or weakly $x$-sparse.
        \end{itemize}

Following \cite{nss7,hjz}, an \emph{$(\ell,w)$-comb} in $G$ is a sequence $((a_i,B_i):i\in[\ell])$ such that $(B_1,\ldots,B_\ell)$ is an $(\ell,w)$-blockade, the vertices $a_1,\ldots,a_\ell$ are distinct and do not belong to $B_1\cup\cdots\cup B_\ell$, and for all distinct $i,j$ the apex $a_i$ is complete to $B_i$ and anticomplete to $B_j$.
A \emph{special vertex} for the comb is a vertex $v$ outside the comb that is complete to $B_1\cup\cdots\cup B_\ell$ and anticomplete to $\{a_1,\ldots,a_\ell\}$. 

\subsection{Results we use}
    The following results are used later. 

\begin{theorem}[R\"odl \cite{rodl}]\label{thm:rodl}
For every graph $H$ and every $\eps\in(0,\frac12)$ there exists $\delta>0$ such that every $H$-free graph $G$ has an $\eps$-restricted induced subgraph with at least $\delta|G|$ vertices.
\end{theorem}

\begin{theorem}[Statement 3.1 of Nguyen--Scott--Seymour \cite{nss5}]\label{thm:sparsepath}
Let $s\ge2$ and $0<z\le 1/(60s)$. Every $z^2$-sparse $P_s$-free graph $G$ with $|G|\ge z^{-2}$ has an anticomplete $(1/z,\lfloor z^2|G|\rfloor)$-blockade.
\end{theorem}
(In \cite{nss5} there is no lower bound on $|G|$, and then the blocks may be empty. We add the hypothesis $|G|\ge z^{-2}$, which gives $\lfloor z^2|G|\rfloor\ge1$, so that the blocks are nonempty.
With this extra hypothesis the statement is a special case of the one in \cite{nss5}.)

\begin{lemma}[Lemma 4.1 of Nguyen--Scott--Seymour \cite{nss7}]\label{lem:cores}
Let $k\ge2$ be an integer, and let $G$ be a graph whose anticonnected components all have fewer than $|G|/k$ vertices. Then $G$ has a complete $(k,|G|/k^2)$-blockade.
\end{lemma}

\begin{lemma}[Lemma 2.10 of Huang--Ju--Zhou \cite{hjz}, contained in Lemma 5.2 of \cite{nss7}]\label{lem:comb}
Let $0<x\le y\le2^{-8}$. For every $y^3$-sparse graph $G$ with $|G|\ge y^{-4}$, one of the following holds:
\begin{itemize}
\item there are disjoint $X,Z\subseteq V(G)$ with $|X|\ge y^4|G|$ and $|Z|\ge(1-4y)|G|$ such that $Z$ is $x$-sparse to $X$;
\item $G$ is $2y^4$-sparse; or
\item for some integer $\ell\in[y^{-1},x^{-2}]$ there is an $(\ell,y^4|G|/\ell^2)$-comb in $G$ with a special vertex.
\end{itemize}
\end{lemma}
    The first two outcomes are directly useful for iterative sparsification. 
    The comb is the main case, and in the proof of the local alternative for $\overline{P_s}$-free graphs it is the only case that uses the forbidden graph. 
    In \cite{nss7,hjz} the comb is turned into a pure blockade.
    For $P_s$ we turn it into a weakly semisparse blockade or a restricted set in Lemma~\ref{lem:combstep}.

\begin{theorem}[Alon--Pach--Solymosi \cite{aps}, see also {\cite[Theorem~1.4]{hjz}}]\label{thm:aps}
If $H_1$ and $H_2$ have the Erd\H{o}s--Hajnal property, then so does every graph obtained from $H_1$ by substituting $H_2$ for a vertex.
\end{theorem}

\begin{theorem}[Chv\'atal \cite{chvatal}]\label{thm:chvatal}
Let $G$ be a graph with a linear order $<$ of $V(G)$ such that there is no induced path $abcd$ of $G$ with $a<b$ and $d<c$. Then $\chi(G)=\omega(G)$.
In particular $G$ has a clique or a stable set of size at least $|G|^{1/2}$.
\end{theorem}
(The second statement follows from $|G|\le\alpha(G)\chi(G)=\alpha(G)\omega(G)$.) A graph with a linear order as in Theorem~\ref{thm:chvatal} is called \emph{perfectly orderable} \cite{chvatal}.

\subsection{The reduction of Huang, Ju and Zhou}\label{subsec:hjz}
The last step of our proof is the reduction of Huang, Ju and Zhou \cite{hjz}, which refines the final part of \cite{nss7}. We recall their definitions, specialised to a single forbidden graph and to a hereditary class.

\begin{definition}[\cite{hjz}]\label{def:gennice}
Let $\cC$ be a hereditary class of graphs. We say $\cC$ is \emph{generalized nice} if there are constants $c_1\ge3$, $c_2\ge8$, $c_3,c_4,c_5,c_8>0$, $c_6\ge1$, $c_7\ge4$ such that for every $G\in\cC$ and every $0<\eps<\frac12$, either
\begin{itemize}
\item $G$ has an $(\eps^{-1},\eps^{c_1}|G|)$-blockade whose blocks are pairwise complete or weakly $\eps^{c_2}$-sparse;
\item $G$ has a clique or stable set of size at least $(\eps^{c_3}|G|)^{c_4}$;
\item $G$ has a complete or anticomplete $(k,|G|/k^{c_5})$-blockade with $k\ge\eps^{-c_6}$; or
\item $G$ has an $\eps^{c_7}$-restricted induced subgraph with at least $\eps^{c_8}|G|$ vertices.
\end{itemize}
\end{definition}
    Our proof will produce only the first two outcomes.
    We retain all four outcomes to match the definition in \cite{hjz}.

\begin{definition}[\cite{hjz}]\label{def:wonderful}
A graph $H$ is \emph{wonderful} if there is a constant $a\ge6$ such that the following holds for every $y\in(0,\frac12)$ and every $\overline{H}$-free graph $G$. If $(B_1,\ldots,B_\ell)$ is an $(\ell,w)$-blockade in $G$ with $\ell\ge y^{-a}$, whose blocks are anticonnected and of the same size, and every two blocks are complete or $y^a$-sparse to each other, then either
\begin{itemize}
\item $G$ has a $y^4$-restricted induced subgraph with at least $w$ vertices; or
\item there is $i\in[\ell]$ such that at most $y|G|$ vertices $v\in V(G)\setminus(B_1\cup\cdots\cup B_\ell)$ satisfy $0<|N(v)\cap B_i|<\frac12|B_i|$.
\end{itemize}
\end{definition}
The graphs in this definition are $\overline H$-free. For $H=P_s$ they are the graphs of the class of $\overline{P_s}$-free graphs, so densities and blockades are in $G$, while the induced copies of $P_s$ that give contradictions are found in $\G$.

A graph $H$ is \emph{leaf-reducible} if it has a vertex $v$ of degree one such that $H\setminus v$ has the Erd\H{o}s--Hajnal property \cite{hjz}. For example, $P_s$ is leaf-reducible whenever $\EH(P_{s-1})$ holds.

\begin{theorem}[A hereditary version of Lemma 1.12 of \cite{hjz}]\label{thm:hjz}
Let $H$ be leaf-reducible and wonderful, and let $\cC$ be a hereditary class of $\overline H$-free graphs. If $\cC$ is generalized nice, then $\cC$ has the Erd\H{o}s--Hajnal property.
\end{theorem}
    
    In the reduction in \cite[Section~3]{hjz}, generalized niceness is applied only to the graph under consideration and its induced subgraphs.
    These graphs remain in $\cC$ by heredity, while the remaining hypotheses follow from the leaf-reducibility and wonderfulness of $H$ and the fact that every graph in $\cC$ is $\overline H$-free.
    When the complement is sparse, \cite[Lemma~2.7]{hjz} applies to the $H$-free complement, without requiring it to belong to $\cC$.
    All constants depend only on $H$ and the constants in Definition~\ref{def:gennice}, so the resulting Erd\H{o}s--Hajnal exponent is uniform over $\cC$.

\section{Wonderfulness of paths}

    In this section we prove that $P_s$ is wonderful whenever $\mathrm{EH}(P_{s-1})$ holds. 
    This is the first of our two new ingredients. 
    Throughout this section, let $s\ge3$ be an integer and put $h=s-1$.

    We begin with some definitions. 
    For a graph $H$ and an integer $c\ge1$, we define $R(H,c)$ recursively.
    Set $R(H,1)=H$.
    For $c\ge2$, let $R(H,c)$ be the graph obtained from $H$ by substituting a copy of $R(H,c-1)$ for every vertex. 
    We also define a graph $H^{(h)}$ as follows.
    First construct pairwise vertex-disjoint graphs $S_h,S_{h-1},\ldots,S_1$ in this order.
    Set $c_h=2^{h-1}$ and let $S_h$ be a copy of $R(H,c_h)$.
    For $i=h-1,h-2,\ldots,1$, once $S_{i+1},\ldots,S_h$ have been constructed, set $c_i=2^{\,i-1+\sum_{j=i+1}^{h}|S_j|}$ and let $S_i$ be a copy of $R(H,c_i)$.
    Let $H^{(h)}$ be obtained from the path $P_h=1\,2\cdots h$ by substituting $S_i$ for vertex $i$ for each $i\in[h]$.
    Thus, for distinct $i,j\in[h]$, the sets $V(S_i)$ and $V(S_j)$ are complete to each other if $|i-j|=1$, and anticomplete otherwise.
    We call $S_1,\ldots,S_h$ the {\em fibres} of $H^{(h)}$.

    We then give two useful properties. 
    \begin{lemma}\label{lem:ramsey}
        Let $H$ be a graph and let $c\ge1$ be an integer.
        Every coloring of $V(R(H,c))$ with at most $c$ colors has a monochromatic induced copy of $H$.
    \end{lemma}

    \begin{proof}
        We prove this lemma by induction on $c$.
        If $c=1$, then $R(H,1)=H$, and the conclusion holds.
        Suppose that $c\ge2$, and let $\varphi$ be a $c$-coloring of $R(H,c)$.
        For each $v\in V(H)$, let $R_v$ be the copy of $R(H,c-1)$ indexed by $v$ in the construction of $R(H,c)$.
        If some $R_v$ uses at most $c-1$ colors, then by the induction hypothesis, $R_v$ contains a monochromatic induced copy of $H$. 
        So we may assume that every $R_v$ uses all $c$ colors.
        For each $v\in V(H)$, choose a vertex $x_v\in V(R_v)$ with $\varphi(x_v)=1$.
        By construction, $\{x_v:v\in V(H)\}$ induces a copy of $H$ whose vertices all have color $1$. 
    \end{proof}

    \begin{lemma}\label{lem:homogenize}
        For every spanning subgraph $M$ of $H^{(h)}$, there are induced copies $H_i$ of $H$ in $S_i$ ($i\in[h]$) such that every pair $(H_i,H_j)$ is complete or anticomplete in $M$. 
    \end{lemma}
    
    \begin{proof}
        For each $i\in[h]$, put $T_i=\bigcup_{j=i+1}^h V(S_j)$.
        We choose induced copies $H_i$ of $H$ in $S_i$ successively for $i=1,\ldots,h$.
        Throughout the construction, we require that $N_M(v)\cap T_j$ is the same for every $v\in V(H_j)$ whenever $H_j$ has been chosen.

        Suppose that $H_1,\ldots,H_{i-1}$ have been chosen.
        For each $j<i$, since $V(S_i)\subseteq T_j$, the sets $N_M(w)\cap V(S_i)$ are equal for all $w\in V(H_j)$.
        Consequently, every vertex $v\in V(S_i)$ is complete or anticomplete to $V(H_j)$ in $M$.
        Let $b_j(v)=1$ if $v$ is complete to $V(H_j)$ in $M$, and let $b_j(v)=0$ otherwise.
        Color each vertex $v\in V(S_i)$ by the tuple $(b_1(v),\ldots,b_{i-1}(v),N_M(v)\cap T_i)$.
        There are at most $2^{i-1}$ possibilities for the first $i-1$ entries and at most $2^{|T_i|}$ possibilities for the last entry.
        Thus at most $2^{i-1+|T_i|}=c_i$ colors occur.
        Since $S_i$ is a copy of $R(H,c_i)$, Lemma~\ref{lem:ramsey} gives a monochromatic induced copy $H_i$ of $H$ in $S_i$.

        For each $j<i$, the value $b_j(v)$ is constant over $v\in V(H_i)$, so $V(H_i)$ and $V(H_j)$ are complete or anticomplete to each other in $M$.
        Moreover, $N_M(v)\cap T_i$ is the same for every $v\in V(H_i)$, preserving the required condition.
        Continuing to $i=h$ gives the desired copies.
    \end{proof}
    \begin{definition}\label{def:lift}
        Let $Q$ be a graph with vertex set $[q]$. 
        A \emph{two-labeled lift} of $Q$ is a graph $L$ with vertex set $\{v,z_1,\ldots,z_q,b_1,\ldots,b_q\}$ (all distinct) such that
        \begin{itemize}
            \item $v$ is adjacent to every $z_i$ and to no $b_i$;
            \item for distinct $i,j$: $z_iz_j$ and $z_ib_j$ are edges if and only if $ij\in E(Q)$;
            \item for distinct $i,j$: if $b_ib_j$ is an edge then $ij\in E(Q)$.
        \end{itemize}
        The adjacency between $z_i$ and $b_i$ is arbitrary, and so is the set of edges $b_ib_j$ with $ij\in E(Q)$. 
        We call $z_1,\ldots,z_q$ the \emph{$z$-vertices}, $b_1,\ldots,b_q$ the \emph{$b$-vertices}, and $i$ the \emph{index} of $z_i$ and of $b_i$. 
    \end{definition}
    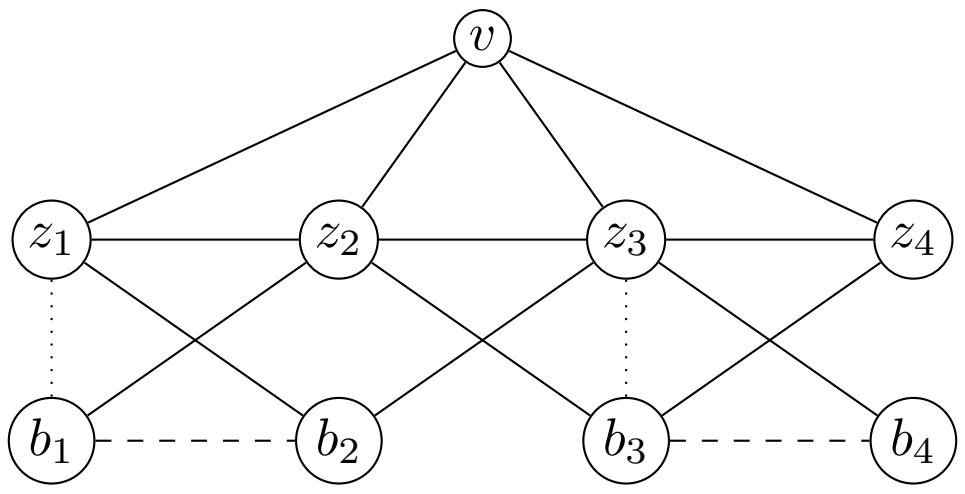
\begin{figure}[ht]
        \centering
        \begin{tikzpicture}[scale=0.95,every node/.style={circle,draw,fill=white,inner sep=1.4pt,font=\small}]
        \node (v) at (3,2.4) {$v$};
        \foreach \i/\x in {1/0,2/2,3/4,4/6} {
            \node (z\i) at (\x,1) {$z_{\i}$};
            \node (b\i) at (\x,-0.4) {$b_{\i}$};
            \draw (v)--(z\i);
        }
        \draw (z1)--(z2) (z2)--(z3) (z3)--(z4);
        \draw (z1)--(b2) (z2)--(b1) (z2)--(b3) (z3)--(b2) (z3)--(b4) (z4)--(b3);
        \draw[dashed] (b1)--(b2) (b3)--(b4);
        \draw[dotted] (z1)--(b1) (z3)--(b3);
    \end{tikzpicture}
    \caption{A two-labeled lift of $Q=P_4$. Solid edges are forced by $Q$. Dashed $b_ib_j$ edges are optional (allowed only on edges of $Q$). Dotted $z_ib_i$ edges are arbitrary.}
    \label{fig:lift}
    \end{figure}
    An \emph{$r$-arm} in a two-labeled lift is an induced path $p_1p_2\cdots p_r$ of the lift such that $p_1$ is a $z$-vertex and $p_2,\ldots,p_r$ are $b$-vertices. 
    Since $v$ is adjacent to $p_1$ and to no $b$-vertex, $vp_1p_2\cdots p_r$ is an induced $P_{r+1}$. 
    The arm is \emph{index-transversal} if its vertices have distinct indices. 
    We then give a graph $Q_h$ such that $\mathrm{EH}(Q_h)$ if $\mathrm{EH}(P_h)$ and prove that every two-labeled lift of $Q_h$ contains an induced $P_{h+1}$. 

    \begin{lemma}\label{lem:forcing}
        For every $h\ge2$ there is a graph $Q_h$ such that 
        \begin{itemize}
            \item[$(1)$] $\mathrm{EH}(Q_h)$ if $\mathrm{EH}(P_h)$; and 
            \item[$(2)$] every two-labeled lift of $Q_h$ contains an induced $P_{h+1}$.
        \end{itemize}
    \end{lemma}
    
    \begin{proof}
        Let $Q_2=P_2$ on $\{1,2\}$. 
        For $2\leq r<h$, let $Q_{r+1}=H^{(h)}$ with $H=Q_r$. 
        We claim that $Q_h$ suffices. 
        Since $Q_h$ is obtained from $P_h$ and $P_2$ by repeated substitution, by \cite{aps}, $\mathrm{EH}(Q_h)$ if $\mathrm{EH}(P_h)$.  
        
        For $2\le r\le h$, say that a graph $Q$ is \emph{$r$-forcing} if every two-labeled lift of $Q$ contains an induced $P_{h+1}$ or an index-transversal $r$-arm. 
        We claim that for $2\le r\le h$, $Q_r$ is $r$-forcing. 
        Since $Q_2=P_2$, in every two-labeled lift of $Q_2$, $z_1b_2$ is an edge. 
        So $z_1b_2$ is an index-transversal $2$-arm. 
        Thus $Q_2$ is $2$-forcing. 
        Now let $2\le r<h$ and suppose $Q_r$ is $r$-forcing. 
        We then prove $Q_{r+1}$ is $(r+1)$-forcing. 

        Suppose to the contrary that $Q_{r+1}$ is not $(r+1)$-forcing. 
        Let $L$ be a two-labeled lift of $Q_{r+1}$ with no induced $P_{h+1}$. 
        Let $M$ be the graph on $V(Q_{r+1})$ in which distinct indices $i,j$ are adjacent if and only if $b_ib_j\in E(L)$.
        So $M$ is a spanning subgraph of $Q_{r+1}$. 
        By Lemma~\ref{lem:homogenize}, there are induced copies $H_1,\ldots,H_h$ of $Q_r$, with $H_i$ in the $i$-th fibre of $Q_{r+1}$, that are pairwise complete or anticomplete in $M$. 
        Note that the vertices of $L$ with indices in $H_i$ (containing both $b$-vertices and $z$-vertices), together with $v$, induce a two-labeled lift of $Q_r$. 
        So each $H_i$ carries an index-transversal $r$-arm. 
        For $i<h$, let $e_i=1$ if $H_i$ is complete to $H_{i+1}$ in $M$ and $e_i=0$ if $H_i$ is anticomplete to $H_{i+1}$ in $M$. 

        \begin{claim}\label{clm:in lem:forcing}
            $e_1=\cdots=e_{r-1}=1$. 
        \end{claim}
        \begin{subproof}[Proof of Claim \ref{clm:in lem:forcing}]
            Suppose $e_i=0$ for some $i\le r-1$. 
            Take an $r$-arm $z',b'_1,\ldots,b'_{r-1}$ with indices in $H_i$, any $b$-vertex $\hat b$ with index in $H_{i+1}$, and a $z$-vertex $\hat z_j$ with index in $H_j$ for each $j=i+2,\ldots,h$. 
            Then $b'_{r-1},\ldots,b'_1,\,z',\,\hat b,\,\hat z_{i+2},\ldots,\hat z_h$ is an induced path of $L$. 
            Its vertices have distinct indices (the arm is index-transversal, and the other vertices lie in other fibres). 
            Recall that two indices in distinct fibres are adjacent in $Q_{r+1}$ if and only if their fibres are consecutive.
            Also, for distinct indices $k,m$, the pair $b_kb_m$ is an edge of $L$ if and only if $km\in E(M)$.
            We check all pairs:
            \begin{itemize}
                \item $b'_{r-1},\ldots,b'_1,z'$ is an induced path, being an arm;
                \item $z'\hat b$ is an edge, since their indices lie in the consecutive fibres $i$ and $i+1$;
                \item $\hat b$ is nonadjacent to every $b'_k$, since $(H_i,H_{i+1})$ is anticomplete in $M$ ($e_i=0$);
                \item $\hat b\hat z_{i+2}$ and $\hat z_j\hat z_{j+1}$ ($i+2\le j<h$) are edges, since their indices lie in consecutive fibres;
                \item every other pair consists of some $\hat z_j$ and a vertex whose index lies in a fibre at distance at least two from fibre $j$, so it is a nonedge.
            \end{itemize}
            The path has $r+h-i\ge h+1$ vertices (as $i\le r-1$), so $L$ contains an induced $P_{h+1}$, a contradiction.
        \end{subproof}
        
        By Claim \ref{clm:in lem:forcing}, choose a $z$-vertex $\hat z$ with index in $H_{r+1}$ (this fibre exists as $r+1\le h$) and a $b$-vertex $\hat b_i$ with index in $H_i$ for $i=1,\ldots,r$. 
        Then $\hat z,\hat b_r,\hat b_{r-1},\ldots,\hat b_1$ is an induced path: consecutive $b$-vertices are adjacent since $e_i=1$, nonconsecutive ones lie in nonconsecutive fibres, and $\hat z$ is adjacent to $\hat b_r$ only. 
        Its vertices have distinct indices, so it is an index-transversal $(r+1)$-arm. 
        This contradicts the choice of $L$, so $Q_{r+1}$ is $(r+1)$-forcing.
        
        Since an $h$-arm together with $v$ is an induced $P_{h+1}$ in $L$, the graph $Q_h$ is then as required. 
        This proves Lemma \ref{lem:forcing}. 
    \end{proof}

    We are now ready to prove $P_s$ is wonderful if $\mathrm{EH}(P_h)$. 
    
    \begin{lemma}\label{lem:wonderful}
        $P_s$ is wonderful if $\mathrm{EH}(P_h)$.
    \end{lemma}
    
    \begin{proof}
        Let $Q_h$ be the graph given by Lemma \ref{lem:forcing} and let $q_h=|Q_h|$, $c$ be the Erd\H{o}s--Hajnal constant of $Q_h$. 
        We claim that $a=\max\{2q_h,2+5/c,6\}$ suffices. 
        Let $y\in (0,1/2)$ and  $G$ be $\overline{P_s}$-free.
        Let $\mathcal{B}=(B_1,\ldots,B_\ell)$ be a $(\ell,w)$-blockade where $\ell$ is an integer at least $y^{-a}$, each block of $\mathcal{B}$ has the same size, and every two blocks are complete or $y^{a}$-sparse.

        Let $I(v)$ be the set of indices $i\in[\ell]$ such that $|N(v)\cap B_i|\in(0,\frac{1}{2}|B_i|)$.
        If $|I(v)|\leq y\ell$ for every vertex $v\in V(G)\setminus V(\mathcal{B})$, then the second bullet of wonderfulness holds by double counting. 
        So we may suppose that there exists $v\in V(G)\setminus V(\mathcal{B})$ with $|I(v)|\geq y\ell$. 
        Without loss of generality, we may assume that $I(v)=[|I(v)|]$. 
        Let $J$ be the graph with $V(J)=I(v)$ and $ij\in E(J)$ if and only if $B_i$ is complete to $B_j$. 

        \begin{claim}\label{clm:pattern-EH}
            The graph $J$ is $\overline{Q_h}$-free.
            Consequently, $J$ has a clique or a stable set of size at least $\lceil|J|^c\rceil$.
        \end{claim}

        \begin{subproof}[Proof of Claim~\ref{clm:pattern-EH}]
            Let $W$ be the common size of the blocks.
            For each $i\in I(v)$, choose $b_i\in N_G(v)\cap B_i$ and let $C_i=B_i\setminus N_G(v)$.
            By the definition of $I(v)$, the vertex $b_i$ exists and $|C_i|>W/2$.

            Suppose that $J[I']\cong\overline{Q_h}$ for some $I'\subseteq I(v)$.
            Identify the vertices of $Q_h$ with the indices in $I'$ so that, for distinct $i,j\in I'$, we have $ij\in E(Q_h)$ if and only if $ij\notin E(J)$.
            Thus, if $ij\in E(Q_h)$, then $B_i$ and $B_j$ are not complete to each other and hence are $y^a$-sparse to each other.
            If $ij\notin E(Q_h)$, then $B_i$ and $B_j$ are complete to each other.

            We choose vertices $z_i\in C_i$ for $i\in I'$ one at a time.
            When choosing $z_i$, exclude the neighbors in $G$ of every $b_j$ with $j\in I'$ and $ij\in E(Q_h)$, and also the neighbors in $G$ of every previously chosen $z_j$ with $ij\in E(Q_h)$.
            Each such vertex has at most $y^a\cdot W$ neighbors in $B_i$, since $B_j$ is $y^a$-sparse to $B_i$.
            There are at most $2(q_h-1)$ such vertices, so at most $2(q_h-1)y^a\cdot W$ vertices of $C_i$ are excluded.
            Since $y<1/2$ and $a\ge2q_h$, we have $y^a<2^{-a}\le2^{-2q_h}$ and hence $2(q_h-1)y^a<2q_h2^{-2q_h}\le1/2$.
            Thus fewer than $W/2$ vertices are excluded, and the choice of $z_i$ is possible because $|C_i|>W/2$.

            Let $L$ be the subgraph of $\overline G$ induced by $\{v\}\cup\{z_i,b_i:i\in I'\}$.
            These vertices are distinct, since the blocks are pairwise disjoint, $v$ lies outside the blockade, and $z_i\notin N_G(v)$ whereas $b_i\in N_G(v)$.
            We verify that $L$ is a two-labeled lift of $Q_h$.
            By construction, $v$ is adjacent in $L$ to every $z_i$ and to no $b_i$.
            For each $ij\in E(Q_h)$, the choices of $z_i$ and $z_j$ ensure that $z_iz_j$, $z_ib_j$, and $z_jb_i$ are edges of $L$.
            For distinct $i,j\in I'$ with $ij\notin E(Q_h)$, the blocks $B_i$ and $B_j$ are complete in $G$, so there are no edges of $L$ between $\{z_i,b_i\}$ and $\{z_j,b_j\}$.
            No further condition is required on the pairs $z_ib_i$, or on the pairs $b_ib_j$ with $ij\in E(Q_h)$.
            Thus $L$ satisfies the definition of a two-labeled lift.
            By Lemma~\ref{lem:forcing}, $L$ contains an induced $P_s$, contrary to $G$ being $\overline{P_s}$-free.
            Therefore $J$ is $\overline{Q_h}$-free.

            It follows that $\overline J$ is $Q_h$-free.
            By the choice of $c$, the graph $\overline J$, and hence also $J$, has a clique or a stable set of size at least $|J|^c$.
        \end{subproof}

        By Claim \ref{clm:pattern-EH}, $J$ has a clique or a stable set $R$ of size $r=\lceil|J|^c\rceil$.
        By taking complement if necessary, we may assume that $R$ is a stable set. Let $S=\bigcup_{i\in R}B_i$. 
        Since each block of $\mathcal{B}$ has the same size, each vertex $v\in B_i\subseteq S$ has at most $r^{-1}|S|$ neighbors in $B_i$ and at most $y^a|S|$ neighbors in $S\setminus B_i$. Since 
        \begin{align*}
            r^{-1}+y^a & \leq (y\ell)^{-c}+y^a  & (r\ge (y\ell)^c)\\
                   &  \leq y^{c(a-1)}+y^a   & (\ell \geq y^{-a}) \\
                   &  \leq y^5+y^5          & (a\ge \frac{5}{c}+2\ge 6) \\
                   &   \leq y^4,            & (y<1/2)
        \end{align*}
        $S$ is $y^4$-restricted. 
        This gives the first outcome of Lemma \ref{lem:wonderful}.
    \end{proof}

\section{Reducing the Erd\H{o}s--Hajnal property to the local alternative}\label{sec:machinery}

    This section turns a local statement about sparse graphs into the Erd\H{o}s--Hajnal property, by the iterative sparsification method of \cite{nss7} (see also \cite{hjz}). 

    \begin{definition}\label{def:LA}
        Let $\cC$ be a hereditary class. We say $\cC$ has the \emph{local alternative with constant $D\ge4$} if for all $0<x\le y\le2^{-D}$, every $y^3$-sparse $F\in\cC$ with $|F|\ge x^{-D}$ satisfies one of the following:
        \begin{itemize}
            \item[\textup{(O1)}] $F$ is $2y^4$-sparse;
            \item[\textup{(O2)}] $F$ has a weakly $x$-semisparse $(k,|F|/k^D)$-blockade for some integer $k\in[1/y,x^{-D}]$;
            \item[\textup{(O3)}] there are disjoint $X,Z\subseteq V(F)$ with $|X|\ge y^D|F|$ and $|Z|\ge(1-4y)|F|$ such that $Z$ is $x$-sparse to $X$;
            \item[\textup{(O4)}] for some $z\in[x^D,y]$ the graph $F$ has a $z^4$-restricted set of at least $z^D|F|$ vertices.
        \end{itemize}
    \end{definition}

    The main result of this section is the following. 
    \begin{theorem}\label{thm:machinery}
        Let $s\ge3$ and assume $\EH(P_{s-1})$. 
        Let $\cC$ be a hereditary class of $\overline{P_s}$-free graphs with the local alternative with constant $D\ge4$. 
        Then $\cC$ has the Erd\H{o}s--Hajnal property.
    \end{theorem}

    Throughout this section, fix $s$, $\cC$ and $D$ as in Theorem~\ref{thm:machinery}. 
    
    Put $c_0=2^{-12},~y_0=\min\bigl(c_0^4/2,~ 1/(60s),~2^{-D}\bigr),~\lambda=6(D+1)$ and $D^*=12D^2+12D+3$.
    Note that $c_0<e^{-8}$ and $D^*=2D\lambda+3$.

    We first derive weakly semisparse blockades from the local alternative.

    \begin{lemma}[Iterative sparsification]\label{lem:bootstrap}
        Let $0<x<y_0$. 
        Every $c_0y_0^3$-sparse $G\in\cC$ with $n=|G|\ge x^{-D^*}$ has a weakly $x$-semisparse $(k,n/k^{D^*})$-blockade for some integer $k$ with $2\le k\le x^{-(2D+1)}$.
    \end{lemma}

    \begin{proof}
        Call $y\in[x^{2D},y_0]$ \emph{admissible}, if $G$ has a $c_0y^3$-sparse induced subgraph $F$ with $|F|\ge y^\lambda n$. 
        If $y$ is admissible, we call the $c_0y^3$-sparse induced subgraph $F$ with $|F|\ge y^\lambda n$ a \emph{witness} for $y$. 
        Every witness is an induced subgraph of the finite graph $G$.
        So there are at most $2^n$ candidate witnesses. 
        For a fixed $F$, every $y$ in 
        $$[\max\{x^{2D},(\frac{\Delta(F)}{c_0|F|})^{1/3}\},\min\{\ y_0,(\frac{|F|}{n})^{1/\lambda}\}]$$ 
        is admissible for $F$.  
        Since there are finitely many $F$, the admissible values form a finite union of closed intervals, hence a closed set, which contains $y_0$. 
        Let $y$ be the smallest admissible value and $F$ be its witness. 

        Suppose first that $y\ge x$. 
        We build disjoint sets $X_0,X_1,\ldots$ and a decreasing sequence $F=F_0\supseteq F_1\supseteq\cdots$ with $|F_j|\ge(1-4y)^j|F|$, as long as $j<K:=\lceil1/y\rceil$.   
        Since $j<1/y$ and $4y\leq 1/2$, $|F_j|\geq (1-4y)^j|F| \ge e^{-8}|F|>c_0|F|$. 
        So $F_j$ is $y^3$-sparse and $|F_j|\ge c_0y^\lambda n\ge x^{\lambda+1}n\ge x^{-D}$, which satisfies the condition of Definition \ref{def:LA}. 
        Apply the local alternative to $F_j$.
        \begin{itemize}
            \item If Definition \ref{def:LA} (O3) holds for $F_j$, then set $X_j=X$ and $F_{j+1}=Z$. 
            Since $|Z|\ge(1-4y)|F_j|$, the bound on $|F_{j+1}|$ is maintained, and we continue. 
            If this happens for all $j<K$, then $(X_0,\ldots,X_{K-1})$ is a blockade in which every later block is $x$-sparse to every earlier one. 
            In particular it is weakly $x$-semisparse. 
            Since $c_0\ge y^2$ and $y\ge1/K$, its width is at least $c_0y^{D+\lambda}n\ge y^{D+\lambda+2}n\ge n/K^{D^*}$ and its length is $K$ with $2\le K\le2/x\le x^{-2}$. 
            
            \item If Definition \ref{def:LA} (O1) holds for $F_j$, then since $2y^4\le c_0^4y^3$, $F_j$ is $c_0(c_0y)^3$-sparse and $|F_j|\ge c_0y^\lambda n\ge(c_0y)^\lambda n$. 
            Since $c_0y\ge c_0x\ge x^2\ge x^{2D}$, the value $c_0y<y$ is admissible, a contradiction.
            
            \item If Definition \ref{def:LA} (O2) holds for $F_j$, then there is a weakly $x$-semisparse $(k,|F_j|/k^D)$-blockade for some integer $k\in[1/y,x^{-D}]$. 
            Since $c_0\ge y\ge1/k$, its width is at least $c_0y^\lambda n/k^D\ge n/k^{D+\lambda+1}\ge n/k^{D^*}$. 
            
            \item If Definition \ref{def:LA} (O4) holds for $F_j$, then let $T$ be the restricted set with $|T|\ge z^D|F_j|\ge c_0z^{\lambda+D}n$, where $z\in [x^D,y]$.
            If $G[T]$ is $z^4$-sparse, put $y'=z^{7/6}$.
            Then since $z\le c_0^2$ and $\lambda=6(D+1)$, $z^4\le c_0(y')^3$ and $|T|\ge z^{\lambda+D+1}n=(y')^\lambda n$.
            Also $x^{2D}\le z^{7/6}<y$.
            So $y'$ is admissible, contradicting minimality.
            
            Thus, $\overline{G}[T]$ is $P_s$-free and $(z^2)^2$-sparse with $z^2\le1/(60s)$ and $z^4|T|\ge z^{\lambda+D+5}n\ge x^{D(\lambda+D+5)}n\ge x^{-1}\ge2$.
            By Theorem~\ref{thm:sparsepath}, there is an anticomplete blockade of $\overline{G}[T]$, that is, a complete blockade of $G[T]$, of length $K'=\lceil z^{-2}\rceil$ and width at least $\lfloor z^4|T|\rfloor\ge z^4|T|/2\ge z^{\lambda+D+6}n\ge n/(K')^{(\lambda+D+6)/2}\geq n/(K')^{D^*}$ with $2\le K'\le2z^{-2}\le x^{-(2D+1)}$.
        \end{itemize}

        So we may assume that $y<x$. 
        Since $y\geq x^{2D}$, $|F|\ge x^{2D\lambda}n$ and $F$ has maximum degree at most $c_0x^3|F|$. 
        Split $F$ into $K=\lceil1/x\rceil$ disjoint sets of size $\lfloor|F|/K\rfloor\ge x|F|/4$. 
        Every vertex has at most $c_0x^3|F|\le x\cdot x|F|/4$ neighbors in each set, so this is an $x$-semisparse blockade of width at least $x^{2D\lambda+3}n\ge n/K^{D^*}$.
    \end{proof}
    
    As in \cite[Lemma~5.5]{nss7} and \cite[Lemma~4.4]{hjz}, we then use R\"odl's theorem to remove the sparsity hypothesis of Lemma~\ref{lem:bootstrap}.

    \begin{lemma}\label{lem:rodlstart}
        There is a constant $D_0\ge1$ such that for every $0<x<y_0$, every $G\in\cC$ with $|G|\ge x^{-D_0}$ has a weakly $x$-semisparse $(k,|G|/k^{D_0})$-blockade for some integer $k$ with $2\le k\le x^{-(2D+1)}$.
    \end{lemma}

    \begin{proof}
        Let $\delta$ be given by Theorem~\ref{thm:rodl} for $H=\overline{P_s}$ and $\eps=c_0y_0^3$.
        Put $L=\lceil\log_2(1/\delta)\rceil$ and $D_0=D^*+L+\lceil\log_2(2(60s)^2)\rceil$. 
        Let $n=|G|\ge x^{-D_0}$, and let $F$ be an $\eps$-restricted induced subgraph with $|F|\ge\delta n\ge x^{-D^*}$ (as $\delta\ge2^{-L}\ge x^L$ and $D_0\ge D^*+L$). 
        If $F$ is $\eps$-sparse, Lemma~\ref{lem:bootstrap} applies to $F\in\cC$ and gives a weakly $x$-semisparse blockade of $F$, hence of $G$, of some length $k$ with $2\le k\le x^{-(2D+1)}$ and width at least
        \[
            |F|/k^{D^*}\ \ge\ \delta n/k^{D^*}\ \ge\ 2^{-L}n/k^{D^*}\ \ge\ n/k^{D^*+L}\ \ge\ n/k^{D_0},
        \]
        using $\delta\ge2^{-L}$ and $k\ge2$.
        So we may assume that $\G[F]$ is $\eps$-sparse.
        Then it is $P_s$-free and $(60s)^{-2}$-sparse since $\eps\le y_0^2$.
        Since $|F|\ge x^{-D^*}\ge(60s)^2$, Theorem~\ref{thm:sparsepath} with $z=1/(60s)$ gives two blocks forming a complete $(2,\delta n/(2(60s)^2))$-blockade in $G$ and $\delta/(2(60s)^2)\ge2^{-D_0}$.
    \end{proof}

    \begin{lemma}[Layout]\label{lem:layout}
        Let $D_0,\gamma\ge1$ and $x_0>0$ and suppose that for every $0<x<x_0$ every $G\in\cC$ with $|G|\ge x^{-D_0}$ has a weakly $x$-semisparse $(k,|G|/k^{D_0})$-blockade for some integer $k$ with $2\le k\le x^{-\gamma}$. 
        Let $q\ge8D_0+8$ be an integer with $2^{-q}<x_0$. 
        Then for every $G\in\cC$ and every $0<\eps<\frac12$, either
        \begin{itemize}
            \item $\homog(G)\ge\eps^{(q+2)D_0}|G|$; or 
            \item $G$ has a $(\lceil1/\eps\rceil,\eps^{(\gamma q+2)D_0}|G|)$-blockade whose blocks are pairwise complete or weakly $\eps^8$-sparse. 
        \end{itemize}
        In particular $\cC$ is generalized nice.
    \end{lemma}

    \begin{proof}
        Let $n=|G|$.
        If $n<\eps^{-(q+2)D_0}$, then $\eps^{(q+2)D_0}|G|<1$ and so the first outcome holds.
        So we may assume $n\geq \eps^{-(q+2)D_0}$. 
        Put $x=\eps^q$, $\ell=\lceil1/\eps\rceil\le2/\eps$ and $\sigma=\eps^{4D_0}$.

        As in \cite{nss7}, a \emph{layout} is a family of disjoint nonempty sets $A_1,\ldots,A_p$ together with a graph $\Pi$ on $[p]$ (the \emph{pattern}) such that $A_i$ is complete to $A_j$ whenever $ij\in E(\Pi)$. The pattern records the pairs of sets known to be complete. 
        A pair of distinct vertices is \emph{undecided} if both lie in the same $A_i$, and \emph{decided} otherwise. 
        In particular, a pair involving a vertex in no $A_i$ is decided, so discarding vertices from the sets $A_i$ never decreases the number of decided pairs.
        A decided pair is \emph{wrong} if it is an edge of $G$ between some $A_i$ and $A_j$ with $ij\notin E(\Pi)$. 

        Choose a layout with the maximum number $p$ of sets subject to
        \[
            |A_i|\ge\sigma n\ (i\in[p]),\qquad \textstyle\sum_i|A_i|^{1/D_0}\ge n^{1/D_0},\qquad\#\text{wrong}\le x\cdot\#\text{decided}.
        \]
        Such a layout exists since the layout with the single set $V(G)$ qualifies.

        If $p\ge\ell$, take $\ell$ of the sets. 
        For $ij\notin E(\Pi)$, every edge between $A_i$ and $A_j$ is a wrong pair. 
        So the pair $(A_i,A_j)$ has at most $x\cdot\#\text{decided}\le xn^2/2$ edges. 
        Hence density is at most $x/(2\sigma^2)=\eps^{q-8D_0}/2\le\eps^8$. 
        If $ij\in E(\Pi)$, then $A_i$ is complete to $A_j$. 
        This gives the second outcome of Lemma \ref{lem:layout} with width $\sigma n$.

        So we may assume that $p<\ell$.
        Then since $\sum_i|A_i|^{1/D_0}\ge n^{1/D_0}$, some set $A$ has $|A|^{1/D_0}\ge n^{1/D_0}/\ell$ and so $|A|\ge(\eps/2)^{D_0}n\ge\eps^{2D_0}n\ge x^{-D_0}$. 
        By hypothesis $G[A]$ has a weakly $x$-semisparse blockade $(A'_1,\ldots,A'_k)$ with $2\le k\le x^{-\gamma}$ and width at least $|A|/k^{D_0}$. 
        If $k\ge\ell$, any $\ell$ of its blocks give the second outcome of Lemma \ref{lem:layout} since $|A|/k^{D_0}\ge\eps^{2D_0}x^{\gamma D_0}n$ and $x\le\eps^8$. 
        So we may assume that $k<\ell$. 
        Replace $A$ by $A'_1,\ldots,A'_k$, joining $A'_i,A'_j$ in the pattern when they are complete, and giving each $A'_i$ the same pattern relations as $A$ to the other sets. 
        This is again a layout.
        Since $k<\ell\le2/\eps$, each new set has at least $|A|/k^{D_0}\ge\eps^{2D_0}(\eps/2)^{D_0}n\ge\sigma n$ vertices and $\sum_i|A'_i|^{1/D_0}\ge k\cdot(|A|/k^{D_0})^{1/D_0}=|A|^{1/D_0}$.
        Old wrong pairs stay wrong or stop being wrong.
        New wrong pairs join two blocks $A_i',A_j'$ forming a weakly $x$-sparse pair, so there are at most $x|A_i'||A_j'|$ of them, while all these pairs were undecided and are now decided.
        So the new layout satisfies the constraints and has more sets, a contradiction.

        The two outcomes are the first two outcomes of Definition~\ref{def:gennice}, with $c_1=(\gamma q+2)D_0\ge3$, $c_2=8$, $c_3=(q+2)D_0$ and $c_4=1$.
        In particular $\cC$ is generalized nice.
    \end{proof}

    Now we are ready to prove Theorem~\ref{thm:machinery}. 
    \begin{proof}[Proof of Theorem~\ref{thm:machinery}]
        Let $D_0$ be given by Lemma~\ref{lem:rodlstart}. 
        Apply Lemma~\ref{lem:layout} with this $D_0$, $\gamma=2D+1$, $x_0=y_0$ and $q=\max(\lceil8D_0+8\rceil,\lceil\log_2(1/y_0)\rceil+1)$, which satisfies $q\ge8D_0+8$ and $2^{-q}<y_0$. 
        So $\cC$ is generalized nice. 
        Since $\EH(P_{s-1})$, $P_s$ is leaf-reducible and wonderful by Lemma \ref{lem:wonderful}. 
        Then, Theorem~\ref{thm:hjz} with $H=P_s$ shows that $\cC$ has the Erd\H{o}s--Hajnal property.
    \end{proof}

\section{Rooted purification}

    This section contains our second new ingredient, a family of statements $\mathsf{R}_r$ ($r\ge5$) that replaces the purity argument of \cite{nss7} for long paths. 
    
    Throughout this section, densities (sparse, complete, restricted, clean) are measured in $G$, while the forbidden induced paths live in $\overline{G}$. We write ``$\overline{G}$-edge'' and ``induced path of $\overline{G}$'' whenever we mean the complement. 

    \subsection{The statement}

        In what follows, $U$ and $Y$ are disjoint sets of vertices, and the vertices of $U$ are called \emph{roots}. For $x>0$, a set $S\subseteq Y$ is \emph{$x$-clean for $U$} if every root is complete to $S$ or has fewer than $x|S|/4$ neighbors in $S$. (The factor $\frac14$ is only for convenience in the cut tree below.)
        
        Let $u\in V(G)\setminus Y$ and $r\ge2$. 
        A \emph{rooted $P_r$ from $u$ into $Y$} is an induced path of $\overline{G}$ with $r$ vertices, one of whose ends is $u$ and whose other $r-1$ vertices lie in $Y$. 
        We say that $u$ \emph{sends no rooted $P_r$ into $Y$} if there is no such path.
        If $u$ sends no rooted $P_r$ into $Y$, then it sends no rooted $P_r$ into any subset of $Y$, and no rooted $P_{r'}$ into $Y$ for any $r'\ge r$ (an induced $P_{r'}$ starting at $u$ begins with an induced $P_r$ starting at $u$). 
        The following observation is used many times.

        \begin{lemma}[Prefixing]\label{lem:prefix}
            Let $r\geq 3$ and $u\notin Y$ send no rooted $P_r$ into $Y$. 
            Let $Z\subseteq Y$ and $p\in Y\setminus Z$ be such that $u$ is complete to $Z$ and nonadjacent to $p$. 
            Then $p$ sends no rooted $P_{r-1}$ into $Z$.
        \end{lemma}
        
        \begin{proof}
            In $\overline{G}$ the vertex $u$ is adjacent to $p$ and to no vertex of $Z$.
            Thus, a rooted $P_{r-1}$ from $p$ into $Z$, together with $u$, would be a rooted $P_r$ from $u$ into $Y$.
        \end{proof}

        \begin{definition}\label{def:Rr}
            Let $r\ge5$ and let $\cC$ be a hereditary class. We say $\cC$ satisfies $\Rs_r$ \emph{with constant $d\ge1$} if the following holds for all $x,y$ with $0<x\le y\le 2^{-d}$. Let $G\in\cC$ and let $U,Y\subseteq V(G)$ be disjoint with $|Y|=n\ge x^{-d}$, such that every root $u\in U$ sends no rooted $P_r$ into $Y$. 
            Then
            \begin{itemize}
                \item[\textup{(R1)}] there is $S\subseteq Y$ with $|S|\ge y^dn$ that is $x$-clean for $U$; or
                \item[\textup{(R2)}] for some $z\in[x^d,y]$ there is a $z^4$-restricted $T\subseteq Y$ with $|T|\ge z^dn$; or
                \item[\textup{(R3)}] for some integer $K\in[1/y,x^{-d}]$ there is an $x$-semisparse $(K,n/K^d)$-blockade in $G[Y]$.
            \end{itemize}
            We say $\cC$ satisfies $\Rs_r$ if it does so with some constant $d$.
        \end{definition}
        Increasing $d$ weakens every outcome and strengthens every hypothesis, so a statement $\Rs_r$ with constant $d$ also holds with every $d'\ge d$.

        We also need the classes in which the ``rooted condition'' holds along a vertex ordering.

        \begin{definition}\label{def:Er}
            For $r\ge2$ and a hereditary class $\cC$, let $\cO_r(\cC)$ be the class of graphs $G\in\cC$ whose vertices can be ordered $v_1,\ldots,v_n$ so that each vertex $v_i$ sends no rooted $P_r$ into $\{v_{i+1},\ldots,v_n\}$. This class is hereditary.
        \end{definition}

        We are now ready to state the main results in this section. 

        \begin{theorem}\label{thm:R5}
            Every hereditary class satisfies $\Rs_5$ with constant $2127$.
        \end{theorem}

        \begin{theorem}\label{thm:Rstep}
            Let $r\ge6$ and let $\cC$ be a hereditary class. Suppose that $\cC$ satisfies $\Rs_{r-1}$ and $\Rs_{\max(r-2,5)}$, and that $\cO_{r-1}(\cC)$ has the Erd\H{o}s--Hajnal property. Then $\cC$ satisfies $\Rs_r$.
        \end{theorem}

    \subsection{Rooted \texorpdfstring{$P_4$}{P4} and modular decomposition}\label{subsec:modules}

        A \emph{module} of a graph $G$ (called a homogeneous set in \cite{hjz} or geschlossen in \cite{gallai}) is a nonempty set $M\subseteq V(G)$ such that every vertex outside $M$ is complete or anticomplete to $M$. 
        A module is \emph{strong} if it is disjoint from or comparable (under inclusion) with every other module. 
        A strong module is prime if it is connected and its complement is connected.
        We use the following standard facts on modular decomposition \cite{gallai}.

        \begin{lemma}\label{lem:standard facts on modular decomposition}\textnormal{\cite{gallai}}
            The following hold for modules in $G$. 
            \begin{itemize}
                \item The strong modules form a tree $T$ under inclusion, whose top node is $V(G)$ and in which the children of each non-singleton node partition it.
                \item Disjoint modules are complete or anticomplete to each other.
                \item Let $N$ be a nontrivial node of $T$. 
                Then one of the following holds. 
                \begin{itemize}
                    \item If $G[N]$ is disconnected, then the children of $N$ are the vertex sets of its components; or 
                    \item if $\overline{G}[N]$ is disconnected, then the children of $N$ are the vertex sets of its anticonnected components; or 
                    \item if $N$ is prime, then every module of $G$ properly contained in $N$ lies in a single child of $N$.
                \end{itemize}
                
                \item Every module is disjoint from or comparable with every strong module.
            \end{itemize}
        \end{lemma}

        \begin{lemma}\label{lem:p4modules}
            Let $Y\subseteq V(G)$ and $u\notin Y$. 
            Then $u$ sends no rooted $P_4$ into $Y$ if and only if every anticonnected component of $G[N(u)\cap Y]$ is a module of $G[Y]$.
        \end{lemma}
\begin{proof}
Let $C$ be (the vertex set of) an anticonnected component of $G[N(u)\cap Y]$, and suppose $C$ is not a module of $G[Y]$, so some $p\in Y\setminus C$ is mixed on $C$. Then $p\notin N(u)$, since a vertex of $N(u)\cap Y$ outside $C$ is complete to $C$. As $\G[C]$ is connected and $p$ is mixed on $C$, there is a $\G$-edge $qq'$ inside $C$ with $p$ nonadjacent to $q$ and adjacent to $q'$ in $G$. Then $u\,p\,q\,q'$ is an induced path of $\G$, a rooted $P_4$ from $u$ into $Y$. Conversely, if $u\,p\,q\,q'$ is a rooted $P_4$ from $u$ into $Y$, then $q,q'\in N(u)$ lie in the same anticonnected component $C$ of $G[N(u)\cap Y]$, and $p$ is mixed on $C$.
\end{proof}

        \begin{lemma}\label{lem:p4purify}
            Let $0<\xi\le2^{-30}$, let $m$ be an integer with $2^{16}\le m\le2\xi^{-1/2}$, and let $Y\subseteq V(G)$ with $|Y|=n\ge16\xi^{-3}$. 
            Let $U\subseteq V(G)\setminus Y$ be such that every root $u\in U$ sends no rooted $P_4$ into $Y$. 
            Then one of the following holds:
            \begin{itemize}
                \item[\textup{(M1)}] $G[Y]$ has a pure $(K,n/K^4)$-blockade for some integer $K\in[m^{1/4},1/\xi]$;
                \item[\textup{(M2)}] $G[Y]$ has a pure $(K,n/(2m))$-blockade for some integer $K$ with $(8K)^4\le m^3\le(9K)^4$;
                \item[\textup{(M3)}] $G[Y]$ has a complete or anticomplete $(\lfloor1/\xi\rfloor,\xi^2n/m)$-blockade;
                \item[\textup{(M4)}] $G[Y]$ has a complete $(\lfloor1/(\xi m)\rfloor,\xi^3n/4)$-blockade; or
                \item[\textup{(M5)}] there is $S\subseteq Y$ with $|S|\ge n/m$ that is $\xi$-clean for $U$.
            \end{itemize}
        \end{lemma}

        \begin{proof}
            Suppose to the contrary that none of the outcomes of Lemma \ref{lem:p4purify} holds. 
            Note that $m^2\xi\le4$ and $\xi m\le2\xi^{1/2}\le2^{-9}$. 
            Put $\nu=\xi^2n/4,~ \tau=n/m,~ K_1=\lfloor1/\xi\rfloor$. 
            Then $\nu\ge4/\xi$ and $\nu\le\xi\tau\le\tau$.

            Now we choose $K_0$. 
            Let $K_0$ be the least integer with $K_0^4\ge m$, and put $a=K_0-1$. 
            Then $a\ge15$ and $a^4<m$. 
            Also $K_0\le K_1$: from $a^4<m$ and $m^2\xi\le4$ we get $a^8\xi<4$, hence $2a\xi<1$ and $K_0=a+1\le2a<1/\xi$.

            By Lemma \ref{lem:standard facts on modular decomposition}, the strong modules of $Y$ form a tree under inclusion. 
            Let $\mathcal L$ be the set of strong modules of $G[Y]$ with at least $\nu$ vertices. 
            It follows that $\mathcal{L}$ contains $Y$. 
            For $X\in\mathcal L$, call a child of $X$ \emph{small} if it is not in $\mathcal L$, and let the \emph{remainder} $S(X)$ be the union of the small children of $X$. 
            Every vertex of $Y$ lies in the remainder of its lowest ancestor in $\mathcal L$, so the remainders partition $Y$. For nested $X'\subsetneq X$ in $\mathcal L$, the remainder $S(X)$ is disjoint from $X'$, and every vertex outside $X'$ is complete or anticomplete to $X'$. 
            By Lemma~\ref{lem:p4modules}, for every $u\in U$, every anticonnected component of $G[N(u)\cap Y]$ is a module of $G[Y]$.

            \begin{claim}\label{clm:An antichain estimate}
                Let $\mathcal{F}$ be a family of disjoint subsets of $Y$, each of size at most $\tau$. 
                If for every integer $K\in [K_0,K_1]$, fewer than $K$ members of $\mathcal{F}$ have size at least $n/K^4$, then $\textstyle\sum_{A\in\mathcal F}|A|\le(K_0-1)\tau+\frac{n}{K_0-1}+\frac{|\mathcal F|\,n}{K_1^4}$.
            \end{claim}
            
            \begin{subproof}[Proof of Claim \ref{clm:An antichain estimate}]
                We split $\mathcal F$ by size. 
                By the counting condition for $K_0$, fewer than $K_0$ members have size at least $n/K_0^4$, contributing at most $(K_0-1)\tau$. 
                For $K_0\le K<K_1$, the counting condition for $K+1$ ($\le K_1$) says that fewer than $K+1$ members have size at least $n/(K+1)^4$. 
                So at most $K$ members have size in $[n/(K+1)^4,n/K^4)$, and they contribute at most $Kn/K^4=n/K^3$. 
                Summing over $K$ and using $K^{-3}\le\frac1{K(K-1)}=\frac1{K-1}-\frac1K$, these members contribute at most $n\sum_{K\ge K_0}\bigl(\frac1{K-1}-\frac1K\bigr)=\frac{n}{K_0-1}$. 
                The remaining members have size less than $n/K_1^4$, and there are at most $|\mathcal F|$ of them.
            \end{subproof}

            Call a strong module of $Y$ {\em heavy} if it has at least $\tau$ vertices. 
            Since $\tau\geq \nu$, a heavy strong module is also an element of $\mathcal{L}$. 
            A heavy node is {\em minimal} if every child of it is not heavy. 
            Minimal heavy nodes are disjoint modules, hence pairwise complete or anticomplete. 
            If there are at least $K_0$ of them we obtain Lemma \ref{lem:p4purify} (M1) with $K=K_0$, since $\tau=n/m\ge n/K_0^4$. 
            So we may assume that there are at most $a$ minimal heavy nodes.
            
            Similarly, we call a non-heavy node of $\mathcal{L}$ {\em maximal} if its parent node is heavy. 
            The maximal non-heavy nodes are disjoint, have sizes in $[\nu,\tau)$, and number at most $n/\nu$. 
            If for some $K\in[K_0,K_1]$ at least $K$ of them have size at least $n/K^4$, we obtain Lemma \ref{lem:p4purify} (M1). 
            The set of maximal non-heavy nodes satisfies the condition of Claim \ref{clm:An antichain estimate}, which gives their total size by 
            $$an/m+n/a+(n/\nu)\,n/K_1^4\le n/a^3+n/15+64\xi^2n\le n/4;$$
            here $an/m\le n/a^3$ as $a^4<m$, $a\ge15$, and $(n/\nu)n/K_1^4\le(4/\xi^2)(2\xi)^4n$ as $K_1\ge1/(2\xi)$.
            All remainders of non-heavy nodes lie in maximal non-heavy nodes, so the remainders of heavy nodes have total size at least $3n/4$. 
            Every heavy node contains a minimal heavy node, and the heavy nodes containing a fixed minimal heavy node form a chain. Hence some chain $\mathcal K$ of heavy nodes has total remainder size at least $3n/(4a)$.

            \begin{claim}\label{clm:A big remainder.}
                For every $X\in \mathcal K$, $|S(X)|<\tau$. 
            \end{claim}
            \begin{subproof}[Proof of Claim \ref{clm:A big remainder.}]
                Suppose to the contrary that some $X\in\mathcal K$ has $|S(X)|\ge\tau$, and let $S=S(X)$. 
                If $G[X]$ is disconnected, then by Lemma~\ref{lem:standard facts on modular decomposition}, the components of $G[S]$ are small children of $X$, each with fewer than $\nu\le\xi|S|\le|S|/K_1$ vertices.
                Applying Lemma~\ref{lem:cores} to $\overline{G}[S]$ gives an anticomplete $(K_1,|S|/K_1^2)$-blockade in $G[S]$.
                Since $|S|/K_1^2\ge\xi^2n/m$, this gives outcome \textup{(M3)} of Lemma~\ref{lem:p4purify}, a contradiction.
                If $\overline{G}[X]$ is disconnected we obtain a complete blockade in the same way. 
                So we may assume that $X$ is prime. 
                Let $K'=\lfloor1/(\xi m)\rfloor\ge2$. 
                If $S$ is $\xi$-clean for $U$, then $|S|\ge\tau=n/m$ gives Lemma \ref{lem:p4purify} (M5), a contradiction.
                So we may assume that some $u\in U$ is not complete to $S$ but has $q\ge\xi|S|/4$ neighbors in $S$. 
                Let $A$ be an anticonnected component of $G[N(u)\cap Y]$ that meets $S$. 
                By Lemma \ref{lem:p4modules}, $A$ is a module.
                Since $X$ is strong, $A\supseteq X$ or $A\subsetneq X$.
                The former would make $u$ complete to $S$. 
                So $A\subsetneq X$, hence $A$ lies in one child of $X$.
                This child meets $S$, so it is small, and $A\subseteq S$ with $|A|<\nu$ by Lemma \ref{lem:standard facts on modular decomposition}.     
                Consequently $N(u)\cap S$ is the union of the anticonnected components of $G[N(u)\cap Y]$ that meet $S$, and these are exactly the anticonnected components of $G[N(u)\cap S]$, each with fewer than $\nu\le q/K'$ vertices (as $q\xi m\ge\xi^2\tau m/4=\nu$). 
                Lemma~\ref{lem:cores} gives a complete $(K',q/K'^2)$-blockade with $q/K'^2\ge\xi\tau(\xi m)^2/4\ge\xi^3n/4$, which is Lemma \ref{lem:p4purify} (M4), a contradiction. 
            \end{subproof}

            By Claim \ref{clm:A big remainder.}, $|S(X)|<\tau$ for all $X\in\mathcal K$. 
            Let $K$ be the largest integer with $(8K)^4\le m^3$.
            Then $K\ge8$ and $m^3\le(9K)^4$.
            Since $(8Ka)^4\le m^3a^4<m^4$, we have $8Ka<m$ and hence $2K\tau<n/(4a)<3n/(4a)$.
            Since $|S(X)|<\tau$ for all $X\in\mathcal K$, we could combine all $S(X)$ into a blockade $(\Gamma_1,\Gamma_2,\ldots,\Gamma_K)$ such that 
            \begin{itemize}
                \item for each $i\leq K-1$, $\tau<|\Gamma_i|<2\tau$; 
                \item for each $X\in \mathcal{K}$, $S(X)\subseteq \Gamma_i$ for some $i\in [K]$; 
                \item for every pair $X_1,X_2\in \mathcal{K}$ with $X_1\supseteq X_2$ such that $S(X_1)\subseteq \Gamma_i$ and $S(X_2)\subseteq \Gamma_j$, $i\leq j$.  
            \end{itemize}
            Since $2K\tau<3n/(4a)$ and $|\Gamma_i|<2\tau$ for every $i\leq K-1$, $|\Gamma_K|> \tau$.
            For every $i\leq K-1$, $\bigcup_{i<j\leq K}\Gamma_j$ is contained in a module $S$ such that $S\cap \Gamma_i=\emptyset$. 
            So every vertex in $\Gamma_i$ is complete or anticomplete to $\bigcup_{i<j\leq K}\Gamma_j$. 
            By keeping the majority kind in each $\Gamma_i$ for $i\leq K-1$, the resulting sets have at least $\tau/2=n/(2m)$ vertices and form a pure blockade of length $K$, which gives Lemma \ref{lem:p4purify} (M2).
            
            This completes the proof of Lemma \ref{lem:p4purify}. 
        \end{proof}
        

    \subsection{A sampling lemma}\label{subsec:sampling}

        For a graph $G$ with a partition $\mathcal{B}=(B_1,B_2,\ldots,B_n)$ where each $B_i$ is nonempty, we call $B_i$ a {\em layer} of $G$. 
        
        \begin{lemma}\label{lem:sampling}
            Let $0<\kappa\le1$ and $g=\max(2,1/\kappa)$. 
            Let $G$ be a graph whose vertex set has size $M$ and is partitioned into $\mathcal{B}=(B_1,B_2,\ldots,B_n)$ such that for every integer $j$ and every set $S$ of size $j$, if $|S\cap B_i|\leq 1$ for every $i\in [n]$, then $G[S]$ has a clique or stable set of size at least $j^\kappa$. 
            If $$0<\eps\le2^{-4(g+10)} \text{ and } \varphi:=\frac{\max\{|B_i|\}}{M}\le\eps^{4g+4},$$ 
            then $G$ has an $\eps$-restricted set of at least $\eps^{2g+2}M$ vertices.
        \end{lemma}

        \begin{proof}
            Suppose to the contrary that Lemma \ref{lem:sampling} does not hold.
            Put
            \[
                q=\lceil\eps^{-2}\rceil,\quad k=\lceil q^g\rceil,\quad\delta=\frac{q}{512k},\quad L=\ln(1/\eps),\quad\ell=\lceil\eps^{-1}\ln(1/\delta)\rceil .
            \]
            Since $q\le2\eps^{-2}$ and $k\le2^{g+1}\eps^{-2g}$, we have $32k^2\varphi\le2^{2g+7}\eps^4\le1$. 
            Since each $B_i$ is nonempty, $\varphi\ge1/M$, and so $M\ge32k^2$ and $\delta M\ge1$.

            \begin{claim}\label{clm:Many cliques or stable sets of size q}
                By taking complement of $G$ if necessary, there are at least $\frac14(M/(2k))^q$ stable sets of size $q$. 
            \end{claim}
            \begin{subproof}[Proof of Claim \ref{clm:Many cliques or stable sets of size q}]
                Choose a set $X$ of $2k$ distinct vertices uniformly at random. 
                The expected number of pairs of $X$ inside one layer is at most $\binom{2k}2\varphi\le2k^2\varphi\le1/16$, so by Markov's inequality, with probability at least $15/16$ no two vertices of $X$ lie in one layer.
                Then any $k$ vertices of $X$ lie in distinct layers, and so contain a clique or stable set of size $k^\kappa\ge q$. 
                So at least half of the $\binom M{2k}$ sets $X$ contain a clique or stable set of size $q$. 
                A fixed $q$-set lies in exactly $\binom{M-q}{2k-q}$ sets $X$, so the number of $q$-sets that are cliques or stable sets is at least
                \[
                \frac12\binom{M}{2k}\Big/\binom{M-q}{2k-q}=\frac12\binom Mq\Big/\binom{2k}{q}\ge\frac12\left(\frac{M}{2k}\right)^q ,
                \]
                where the identity is $\binom M{2k}\binom{2k}q=\binom Mq\binom{M-q}{2k-q}$ and the inequality holds termwise in $\binom Mq/\binom{2k}q=\prod_{i<q}\frac{M-i}{2k-i}$. 
                By taking complement of $G$ if necessary, at least $\frac14(M/(2k))^q$ of them are stable sets. 
                This proves Claim \ref{clm:Many cliques or stable sets of size q}. 
            \end{subproof}

            \begin{claim}\label{clm:Few stable q-sets}
                If $\ell \leq q$, then the number $I_q$ of stable $q$-sets satisfies $I_q\le\binom{M}{\ell}\binom{N_0}{q-\ell}$, where $N_0=\lceil\delta M\rceil\le2\delta M$. 
            \end{claim}
            \begin{subproof}[Proof of Claim \ref{clm:Few stable q-sets}]
                Note that no set of at least $\delta M$ vertices is $\eps$-sparse in $G$ since $\delta=q/(512k)\geq2^{-(g+9)}\eps^{2g-2}\geq\eps^{2g+2}$, where the last inequality follows from $\eps^4\leq2^{-16(g+10)}\leq2^{-(g+9)}$.  
                Fix a linear order of $V(G)$ and a stable $q$-set $I$. 
                Run the following procedure. 
                Start with $X=V(G)$ and $W=\emptyset$. 
                While $|X|>N_0$, let $p$ be the first vertex of maximum degree in $G[X]$. 
                If $p\notin I$, then delete $p$ from $X$ while if $p\in I$, put $p$ into $W$ and delete $p$ and all its neighbors from $X$. 
                Note that a vertex $v$ of $I$ is deleted if and only if $v\in W$. 
                So, in the end, $I=W\cup (I\cap X)$, where $X$ has at most $N_0$ vertices. 
                Since $G[X]$ is not $\eps$-sparse while $|X|>N_0$, each step when $p\in I$, we delete more than $\eps|X|$ vertices. 
                Since $(1-\eps)^\ell M\le\delta M$, at most $\ell$ vertices of $I$ are put into $W$.  
                Given $W$, we can rerun the procedure without knowing $I$, since at each step, $p\in I$ if and only if $p\in W$. 
                So $W$ determines the final $X$. 
                Pad $W$ to exactly $\ell$ vertices by adding vertices of $I\cap X$. 
                The added vertices lie in the final $X$, so they are never chosen as $p$, hence rerunning the procedure with the padded set in place of $W$ gives the same run and the same final $X$. 
                Then $I$ is determined by an $\ell$-subset of $V(G)$ together with a $(q-\ell)$-subset of $X$. 
                This proves Claim \ref{clm:Few stable q-sets}.
            \end{subproof}

            Put $\Lambda=\ln(1/\delta)$ and $\Lambda'=\ln(1/(2\delta))$.
            Since $\delta=q/(512k)$, $k\leq2q^g$, and $q\leq2\eps^{-2}$, we have $\Lambda\leq(g+9)\ln2+2(g-1)L\leq2gL$, where the last inequality follows from $L\geq4(g+10)\ln2$.

            We first show that $\ell<q$ and $(\ell/q)\Lambda'<1$.
            Since $g\geq2$ and $L\geq4(g+10)\ln2$, we have $L>24$ and $2g<L$.
            Using $\ell/q\leq\eps\Lambda+\eps^2$ and $\eps=e^{-L}$, we obtain
            \[
                2gL\frac{\ell}{q}\leq4g^2\eps L^2+2g\eps^2L\leq2L^4e^{-L}<2\cdot24^4e^{-24}<1,
            \]
            where we used that $t^4e^{-t}$ is decreasing for $t\geq24$.
            As $2gL>1$ and $\Lambda'\leq2gL$, it follows that $\ell<q$ and $(\ell/q)\Lambda'<1$.
            Now by Claim \ref{clm:Few stable q-sets}, the number of stable $q$-sets satisfies 
            \[
                I_q\leq\binom M\ell\binom{N_0}{q-\ell}\leq\binom q\ell\frac{M^\ell N_0^{q-\ell}}{q!}\leq\left(\frac{2e}{q}\right)^qM^\ell(2\delta M)^{q-\ell},
            \]
            where we used $\binom q\ell\leq2^q$, $q!\geq(q/e)^q$, and $N_0\leq2\delta M$.
            Comparing this with $I_q\geq\frac14(M/(2k))^q$ by Claim \ref{clm:Many cliques or stable sets of size q} and taking $q$-th roots yields $(1/(2\delta))^{1-\ell/q}\leq4e(k/q)4^{1/q}\leq8ek/q$, since $q\geq2$.
            Then by taking logarithms and using $\Lambda'=\ln(256k/q)$, we obtain $(1-\ell/q)\Lambda'\leq\ln(8ek/q)=\Lambda'-\ln(32/e)$.
            Hence $(\ell/q)\Lambda'\geq\ln(32/e)>1$, a contradiction.
        \end{proof}

    \subsection{The cut tree}\label{subsec:tree}
        In this subsection, fix $r\ge5$, a graph $G$ and disjoint $U,Y\subseteq V(G)$ with $|Y|=n$ such that every root $u\in U$ sends no rooted $P_r$ into $Y$. 
        Set $x\in(0,1)$, an integer $t\ge2^{16}$ and let $\mu=n/t^6$.

        The \emph{cut tree} is a tree whose nodes are subsets of $Y$ with the children of each node ordered such that the following hold. 
        \begin{itemize}
            \item Its top node is $Y$.
            \item A node $S$ has children if and only if $|S|\ge\mu$ and $S$ is not $x$-clean for $U$. 
            \item The children of a node $S$ are the vertex sets of the anticonnected components of $G[N(u_S)\cap S]$ (in any order), followed by $S\setminus N(u_S)$, where $u_S$ is a root in $U$ such that $x|S|/4\le|N(u_S)\cap S|<|S|$.
        \end{itemize}
        A node $S$ is \emph{large} if $|S|\ge\mu$ and \emph{small} otherwise. 
        We call the vertex sets of the anticonnected components of $G[N(u_S)\cap S]$ the \emph{neighbor children} of $S$ and $S\setminus N(u_S)$ the \emph{non-neighbor child} of $S$. 
        It follows that distinct neighbor children of $S$ are complete to each other and $u_S$ is complete to every neighbor child of $S$ but anticomplete to its non-neighbor child.
        Since children are proper subsets, the tree is finite. 
        Any two nodes are nested or disjoint and the left-to-right order of the tree induces a linear order on any family of pairwise disjoint nodes. 

        \begin{definition}\label{def:frontier}
            A \emph{rooted frontier of order $r$} in $Y$ is a blockade $(Z_1,\ldots,Z_m)$ in $G[Y]$ with nonempty blocks, together with sets $R_1,\ldots,R_m\subseteq V(G)\setminus Y$ of at most $m$ roots each and a graph $J$ with vertex set $[m]$, such that:
            \begin{itemize}
                \item[\textup{(F1)}] every vertex of $R_1\cup\cdots\cup R_m$ sends no rooted $P_r$ into $Y$, and $R_j$ is complete to $Z_j$ for each $j$; and 
            
                \item[\textup{(F2)}] for $i<j$, every vertex of $Z_j$ sends no rooted $P_{r-1}$ into $Z_i$; and 
            
                \item[\textup{(F3)}] if $i<j$ and $ij\notin E(J)$,  then $Z_i$ is complete to $Z_j$; and 
                
                \item[\textup{(F4)}] if $j,\ell$ with $j<\ell$ both lie on a path $j=i_1-i_2-\ldots-i_p=\ell$ in $J$ such that $i_1<i_2<\cdots<i_p$, then some vertex of $R_j$ is anticomplete to $Z_\ell$.
            \end{itemize}
        \end{definition}

        \begin{lemma}\label{lem:tree}
            One of the following holds:
            \begin{itemize}
                \item[\textup{(T1)}] some $S\subseteq Y$ with $|S|\ge\mu$ is $x$-clean for $U$;
                \item[\textup{(T2)}] $G[Y]$ has a complete $(t,n/t^7)$-blockade;
                \item[\textup{(T3)}] there is a rooted frontier of order $r$ in $Y$ with $t$ blocks, each with at least $n/t^6$ vertices;
                \item[\textup{(T4)}] for some integer $m$ with $t\le m\le4t^6/x$, there is a rooted frontier of order $r$ in $Y$ with $m$ blocks, each with at least $n/m^3$ vertices;
                \item[\textup{(T5)}] there are disjoint nonempty sets $L_1,\ldots,L_p\subseteq Y$, each with fewer than $n/t^6$ vertices and with total size at least $3n/(8t)$ such that for every set $X\subseteq L_1\cup \ldots\cup L_p$ with $|X\cap L_i|\leq 1$ for each $i$, the vertex in $X\cap L_i$ (if it exists) sends no rooted $P_{r-1}$ into $X\cap (L_{i+1}\cup \ldots \cup L_p)$. 
                In particular, if $G$ belongs to a hereditary class $\cC$, then $G[X]\in\cO_{r-1}(\cC)$.
            \end{itemize}
        \end{lemma}


        \begin{proof}
            Suppose to the contrary that none of the outcomes of Lemma \ref{lem:tree} holds. 
            Build the cut tree of $Y$. 
            If some large node is $x$-clean for $U$, then we have Lemma \ref{lem:tree} (T1). 
            So we may assume that every node without children is small. 
            We use the following terms.
            \begin{itemize}
                \item A \emph{piece} is a small node whose parent is large. Since every node without children is small, the pieces partition $Y$.
                
                \item The \emph{retained tree} is the subtree formed by the large nodes.
                
                \item A large node is \emph{exceptional} if it has no large child or at least two large children, and \emph{unary} if it has exactly one large child.
                
                \item Let $S$ be a unary node. 
                If its large child is the non-neighbor child, the \emph{positive layer} of $S$ is $N(u_S)\cap S$, the union of its neighbor children. 
                If its large child is a neighbor child, the \emph{negative layer} of $S$ is its non-neighbor child $S\setminus N(u_S)$, and the \emph{fringe} of $S$ is the union of its other small neighbor children.
            \end{itemize}

            \begin{claim}\label{clm:no t partwise disjoint large nodes}
                There are at most $t-1$ pairwise disjoint large nodes. 
            \end{claim}
            \begin{subproof}[Proof of Claim \ref{clm:no t partwise disjoint large nodes}]
                Suppose to the contrary that there are $t$ pairwise disjoint large nodes $Z_1,\ldots,Z_t$, listed in the left-to-right order of the cut tree. 
                In this case, we claim that Lemma \ref{lem:tree} (T3) holds, that is, there is a rooted frontier of order $r$ in $Y$ with $t$ blocks, each with at least $n/t^6$ vertices. 
                For $j<\ell$, let $S$ be the lowest common ancestor of $Z_j$ and $Z_\ell$. 
                Since the non-neighbor child comes last, either $Z_j,Z_\ell$ lie in distinct neighbor children of $S$, and then they are complete to each other, or $Z_j$ lies in a neighbor child, $Z_\ell$ lies in the non-neighbor child and $u_S$ is complete to $Z_j$ but anticomplete to $Z_\ell$. 
                In the second case put $u_S$ into $R_j$ and $j\ell$ into $E(J)$. 
                Then $|R_j|\le t-1$ and (F1) and (F3) hold. 
                For (F2), in the first case a vertex of $Z_\ell$ has no $\G$-neighbor in $Z_j$, and in the second case we use Lemma~\ref{lem:prefix} with $u=u_S$. 
                For (F4), let $i_1<\cdots<i_p$ be an increasing path of $J$, $j=i_1$, $\ell=i_p$ and $S$ be the lowest common ancestor of $Z_j,Z_\ell$. 
                Note that for any block $Z_i\cap S\neq \emptyset$, $Z_i$ lies inside a child of $S$. 
                Since the blocks are listed in the left-to-right order of the cut tree, the indices of the blocks inside $S$ form an interval. 
            
                Suppose $Z_j$ and $Z_\ell$ were in distinct neighbor children $C\ne C'$ of $S$.  
                Then every index between $j$ and $\ell$ belongs to a child lying between $C$ and $C'$ in the order, hence to a neighbor child of $S$.
                In particular all of $i_1,i_{2},\ldots,i_p$ lie in neighbor children of $S$ with $i_1$ in $C$ and $i_p$ in $C'$. 
                Taking the least $c\ge 1$ with $Z_{i_{c+1}}\notin C$, then the blocks $Z_{i_c}$ and $Z_{i_{c+1}}$ lie in distinct neighbor children of their lowest common ancestor $S$.  
                So $i_ci_{c+1}\notin E(J)$ by construction, contradicting that the path uses edges of $J$.
            
                Since the non-neighbor child comes last, we may therefore assume that $Z_j$ is in a neighbor child and $Z_\ell$ in the non-neighbor child, and the root recorded in $R_j$ is anticomplete to $Z_\ell$. This gives (F4) and hence gives Lemma \ref{lem:tree} (T3). 
            \end{subproof}

            Nodes of the retained tree without large children are pairwise disjoint. 
            Hence we may assume that the retained tree has fewer than $t$ nodes without large children by Claim \ref{clm:no t partwise disjoint large nodes}. 
            In a finite rooted tree the number of nodes with at least two children is less than the number of leaves, so the retained tree also has fewer than $t$ nodes with at least two large children. 
            Thus, there are fewer than $2t$ exceptional nodes, and every other large node is unary.

            Suppose some large node has a set $D$ of small neighbor children of total size at least $t\mu$. 
            The anticonnected components of $G[D]$ are these children, each smaller than $\mu\le|D|/t$. 
            By Lemma~\ref{lem:cores}, there is a complete $(t,|D|/t^2)$-blockade with $|D|/t^2\ge\mu/t=n/t^7$, which gives Lemma \ref{lem:tree} (T2). 
            So we may assume that the pieces that are children of any one large node have total size less than $(t+1)\mu$. 
            Then the exceptional nodes have pieces of total size less than $2t(t+1)\mu$, and so the unary nodes have pieces of total size at least $n-2t(t+1)\mu$. 
            Since there are at most $t$ top-to-bottom paths of the retained tree, one such path $\mathcal P$ carries pieces of unary nodes of total size at least $n/t-2(t+1)\mu$.    

            By the definition of fringe sets, fringe sets of different unary nodes of $\mathcal P$ are complete to each other. 
            If the fringe sets of unary nodes of $\mathcal P$ have total size at least $t\mu$, then by Lemma~\ref{lem:cores} again, Lemma \ref{lem:tree} (T2) holds. 
            Since $t\geq 2$ and $\mu=n/t^6$, $n/t-(3t+2)\mu \geq n/t-4t\mu=n/t-4n/(t^5)\geq n/t-n/(4t)=3n/(4t)$. 
            So we may assume that the positive and negative layers of unary nodes of $\mathcal P$ have total size at least $3n/(4t)$. 
            Thus, either the positive layers or the negative layers have total size at least $3n/(8t)$.     
            
            Suppose first that the negative layers $L_1,\ldots,L_p$ (in the order of $\mathcal P$) have total size at least $3n/(8t)$. 
            So $|L_i|<\mu$ for each $i\in [p]$. 
            Let $X$ be a set such that $|X\cap L_i|\leq 1$ and let $p_i\in X\cap L_i$ if $p_i$ exists. 
            The root $u_i$ of the node of $L_i$ is anticomplete to $L_i$ and complete to every later layer. 
            Since $u_i$ sends no rooted $P_r$ into $Y$, by Lemma~\ref{lem:prefix}, $p_i$ sends no rooted $P_{r-1}$ into the later vertices of $X\cap (L_{i+1}\cup \ldots\cup L_p)$. 
            This gives Lemma \ref{lem:tree} (T5), a contradiction. 

            So we may assume that the positive layers have total size at least $3n/(8t)$. 
            Since they are disjoint and each has at least $x\mu/4$ vertices, there are at most $4t^6/x$ of them. 
            We claim that for some integer $m\ge t$, at least $m$ positive layers have at least $n/m^3$ vertices each. 
            Otherwise, listing their sizes as $a_1\ge a_2\ge\cdots$ (with $a_i=0$ beyond the number of positive layers), we have $a_m<n/m^3$ for all $m\ge t$ and $a_i<t\mu=n/t^5$ for $i<t$ since every positive layer has fewer than $t\mu$ vertices. 
            So the total size is
            \[
                \sum_{i\ge1}a_i<\sum_{i<t}\frac{n}{t^5}+\sum_{m\ge t}\frac{n}{m^3}\le\frac{n}{t^4}+n\int_{t-1}^\infty\frac{dz}{z^3}=\frac{n}{t^4}+\frac{n}{2(t-1)^2}<\frac{n}{8t},
            \]
            the last inequality since $t\ge2^{16}$. 
            This contradicts the total size being at least $3n/(8t)$. 
            Keep $m$ such positive layers $Z_1,\ldots,Z_m$ in the order of $\mathcal P$ and let $R_j$ consist of the root of the node of $Z_j$, $J$ be the complete graph. 
            Then the root of $Z_j$ is complete to $Z_j$ and anticomplete to all later positive layers, which gives (F1), (F3) and (F4), and (F2) follows from Lemma~\ref{lem:prefix}. 
            This gives Lemma \ref{lem:tree} (T4), a contradiction. 
            
            This completes the proof of Lemma \ref{lem:tree}.
        \end{proof}

    \subsection{Frontiers of order five}\label{subsec:frontier5}

        \begin{definition}\label{def:FP}
            Let $r\ge5$, $c\ge16$ and $\cC$ be a hereditary class. 
            We say $\cC$ has the \emph{frontier property of order $r$ with constant $c$} if the following holds for every $G\in\cC$, every $Y\subseteq V(G)$, every $n>0$, every integer $m\ge2^c$, every $\rho\in(0,1]$ and every $b\in\{3,6\}$. 
            If $(Z_1,\ldots,Z_m)$ is a rooted frontier of order $r$ in $Y$ with $|Z_i|\ge(m/\rho)^c$ and $|Z_i|\ge n/m^b$ for every $i$, then either
            \begin{itemize}
                \item for some $z\in[(\rho/m)^c,1/m]$, there is a $z^4$-restricted $T\subseteq Z_1\cup\cdots\cup Z_m$ with $|T|\ge z^{c+b}n$; or
                \item for some integer $K\in[m^{1/4},(m/\rho)^c]$, there is a $\rho$-semisparse $(K,n/K^{c+4b})$-blockade in $G[Z_1\cup\cdots\cup Z_m]$.
            \end{itemize}
        \end{definition}

        Let $G$ be a graph with blockade $\mathcal{W}=(W_1,\ldots,W_m)$ and $\beta>0$ such that for every vertex $w\in W_j$ and every $i<j$, $w$ is either complete to $W_i$ or has fewer than $\beta |W_i|$ neighbors in $W_i$. 
        We say that $w$ is {\em complete (sparse) towards} $i$.        
        The {\em ambiguity graph} of $\mathcal{W}$ is the graph $\mathcal{A}$ on $[m]$ in which $i<j$ are adjacent if $W_j$ contains both a vertex complete towards $i$ and a vertex sparse towards $i$. 
        If $\mathcal A$ has no edges, then each $W_j$ is complete or $\beta$-sparse to each earlier $W_i$, and so $(W_1,\ldots,W_m)$ is already $\beta$-semisparse. 
        An \emph{increasing path} of $\mathcal A$ is a sequence $i_1<\cdots<i_p$ with $i_ai_{a+1}\in E(\mathcal A)$ for all $a<p$.

        \begin{lemma}[Monotone paths]\label{lem:monopath}
            
            Let $m\geq 4, k'=\lfloor\sqrt m\rfloor$, $G$ be a graph with blockade $\mathcal{W}=(W_1,\ldots,W_m)$ and $\beta>0$ such that for every vertex $w\in W_j$ and every $i<j$, $w$ is either complete or sparse towards $i$. 
            Let $\mathcal A$ be the ambiguity graph of $\mathcal{W}$.
            Suppose that for every increasing path $i_1<\cdots<i_p$ of $\mathcal A$ and all $a<b-1$, no vertex of $W_{i_b}$ is sparse towards $i_a$ and complete towards $i_{a+1}$. 
            Then there is a $k'\beta$-semisparse blockade of length $k'$ and width at least $\min_{i\in[m]}\{|W_i|/k'\}$, each of whose blocks is a subset of some $W_i$.
        \end{lemma}
        
        \begin{proof}
            By hypothesis, along an increasing path $i_1<\cdots<i_p$, a vertex of $W_{i_b}$ that is sparse towards some $i_a$ is also sparse towards $i_{a+1},\ldots,i_{b-1}$. So its types towards $i_1,\ldots,i_{b-1}$ are complete up to some point and sparse afterwards. 
            Call the number of complete ones its \emph{complete prefix}. 
            For example, a vertex of $W_{i_5}$ whose types towards $i_1,i_2,i_3,i_4$ are complete, complete, sparse, sparse has complete prefix $2$.
            
            Suppose first that $\mathcal A$ has an increasing path $i_1<\cdots<i_{k'}$.    
            Partition each $W_{i_b}$ according to the complete prefix (at most $k'$ values) and keep a largest part $W'_b$, so $|W'_b|\ge|W_{i_b}|/k'$. 
            For every $a<b$, all vertices of $W'_b$ have the same type towards $i_a$. 
            If it is complete, $W'_b$ is complete to $W'_a\subseteq W_{i_a}$ and if it is sparse, every vertex of $W'_b$ has fewer than $\beta|W_{i_a}|\le k'\beta|W'_a|$ neighbors in $W'_a$. 
            Then $(W'_1,\ldots,W'_{k'})$ is the required blockade.

            So we may assume that $\mathcal{A}$ has no increasing path with at least $k'$ vertices. 
            Label each $i\in[m]$ by the number of vertices of a longest increasing path of $\mathcal A$ ending at $i$. 
            So there are at most $k'-1$ labels. 
            Then some label class has at least $m/(k'-1)>k'$ indices. 
            If two indices $i<j$ in this label class are adjacent in $\mathcal A$, then the longest increasing path ending at $j$ is longer than the path ending at $i$, a contradiction. 
            So for every $i<j$ in this class, all vertices of $W_j$ have the same type towards $i$, so $W_j$ is complete or $\beta$-sparse to $W_i$. 
            Hence the corresponding sets $W_i$ form a $\beta$-semisparse blockade of length at least $k'$.
        \end{proof}

        \begin{lemma}\label{lem:frontier5}
            Every hereditary class has the frontier property of order $5$ with constant $20$.
        \end{lemma}

        \begin{proof}
            Let $(Z_1,\ldots,Z_m)$ with $R_j$, $J$ be a rooted frontier of order $5$ with $c=20$ as in Definitions~\ref{def:frontier} and \ref{def:FP}. 
            Put
            \[
                \xi_0=\frac{\rho}{64m^5},\qquad\theta=\frac{\xi_0}4,\qquad \eta=\frac1{16m^3},\qquad \beta=2m^2\theta=\frac{\rho}{128m^3}.
            \]
            Since $m\ge2^{20}$ and $0<\rho\le1$, we have $0<\xi_0<2^{-30}$ and $2^{16}\le m\le2\xi_0^{-1/2}$.
            Moreover, $|Z_i|\ge(m/\rho)^{20}\ge16\xi_0^{-3}$ for every $i\in[m]$.
            By Definition~\ref{def:frontier} (F2), every vertex of $\bigcup_{j>i}Z_j$ sends no rooted $P_4$ into $Z_i$, and this root set is disjoint from $Z_i$.
            Thus the conditions of Lemma~\ref{lem:p4purify} hold for each $Z_i$ with parameters $(\xi_0,m)$ and root set $\bigcup_{j>i}Z_j$.

            In the rest of the proof, whenever an application of Lemma~\ref{lem:p4purify} or Lemma~\ref{lem:cores} gives a blockade instead, we stop at once (an \emph{early outcome}). 
            Such a blockade is pure (Lemma \ref{lem:p4purify} (M1)--(M3)) or complete (\ref{lem:p4purify} (M4) and Lemma~\ref{lem:cores}), hence $\rho$-semisparse. 
            We show that its length and width satisfy the second alternative of Definition~\ref{def:FP} in Table~\ref{tab:frontier5}. 
            So each early outcome is already the conclusion of the lemma, and nothing is carried to the end.

            We proceed as follows.

            \emph{Round 1.} By Lemma~\ref{lem:p4purify} applied to $Z_i$, there is $Z_i^{(1)}\subseteq Z_i$ with $|Z_i^{(1)}|\ge|Z_i|/m$ such that every vertex of $\bigcup_{j>i}Z_j$ is complete to $Z_i^{(1)}$ or has fewer than $\theta|Z_i^{(1)}|$ neighbors in $Z_i^{(1)}$. Accordingly every vertex $p\in Z_j^{(1)}$ has a \emph{type} towards each $i<j$: \emph{complete} or \emph{sparse}. 
            All sets below are subsets of the $Z_i^{(1)}$, and a vertex complete towards $i$ remains complete to every subset of $Z_i^{(1)}$.

            \emph{Round 2.} Call a pair $i<j$ \emph{good} if some root $u\in R_j$ has at least $\eta|Z_i^{(1)}|$ nonneighbors in $Z_i^{(1)}$, and \emph{bad} otherwise. 
            For a good pair fix such a root $u$, and let $S_{ij}$ be the set of vertices of $Z_j^{(1)}$ that are sparse towards $i$. 
            There are at most $\theta|Z_i^{(1)}||S_{ij}|$ edges between $Z_i^{(1)}$ and $S_{ij}$. 
            So some nonneighbor $p_{ij}$ of $u$ in $Z_i^{(1)}$ has at most $(\theta|Z_i^{(1)}||S_{ij}|)/(\eta|Z_i^{(1)}|)\leq\theta|Z_j^{(1)}|/\eta$ neighbors in $S_{ij}$. By Definition~\ref{def:frontier} (F1) and Lemma~\ref{lem:prefix} (with $Z=Z_j^{(1)}$), $p_{ij}$ sends no rooted $P_4$ into $Z_j^{(1)}$. Apply Lemma~\ref{lem:p4purify} to $Z_j^{(1)}$ (which has at least $16m^6$ vertices) with parameters $(m^{-2},m)$ and the representatives $p_{ij}$ ($i<j$ good) as roots. 
            Outside early outcomes we get $D_j\subseteq Z_j^{(1)}$ with $|D_j|\ge|Z_j^{(1)}|/m$ on which each $p_{ij}$ is complete or has fewer than $|D_j|/(4m^2)$ neighbors. 
            
                        After choosing all the sets $D_j$, we make the following deletions from each $D_j$, simultaneously for all pairs.
            \begin{itemize}
                \item For each good pair $i<j$, if $p_{ij}$ is complete to $D_j$, delete $D_j\cap S_{ij}$, which has at most $\theta|Z_j^{(1)}|/\eta\le\theta m|D_j|/\eta$ vertices.
                Otherwise, delete the vertices of $D_j$ that are complete towards $i$.
                All these vertices are adjacent to $p_{ij}$, so there are fewer than $|D_j|/(4m^2)$ of them.

                \item For each bad pair $j<\ell$, delete every vertex of $D_j$ that is nonadjacent to at least one root in $R_\ell$.
            \end{itemize}
            For fixed $j$, the deletions in the second item remove fewer than $m^2\eta|Z_j^{(1)}|\le|D_j|/16$ vertices in total.
            Since $\theta m^2/\eta=\rho/16$, the deletions of the first kind remove at most $m(1/(4m^2)+\theta m/\eta)=1/(4m)+\rho/16\le1/8$ of $D_j$. 
            Let $Z_j^{(2)}$ be what remains of $D_j$. 
            So $|Z_j^{(2)}|\ge|Z_j^{(1)}|/(2m)$. 
            Now for a good pair $i<j$, all vertices of $Z_j^{(2)}$ have the same type towards $i$ and for a bad pair $i<j$ every root in $R_j$ is complete to $Z_i^{(2)}$.

            \emph{Round 3.} If for some $i$, every anticonnected component of $G[Z_i^{(2)}]$ has fewer than $|Z_i^{(2)}|/m$ vertices, Lemma~\ref{lem:cores} gives a complete $(m,|Z_i^{(2)}|/m^2)$-blockade (an early outcome). 
            So there is an anticonnected component of $G[Z_i^{(2)}]$, denoted by $Z_i^{(3)}$, with at least $|Z_i^{(2)}|/m\ge|Z_i|/(2m^3)$ vertices. 
            A vertex sparse towards $i$ has fewer than $\theta|Z_i^{(1)}|\le\beta|Z_i^{(3)}|$ neighbors in $Z_i^{(3)}$. So the types of round~1 satisfy the setting of Lemma~\ref{lem:monopath} for the sets $W_i=Z_i^{(3)}$ with this $\beta$. 
            Let $\mathcal A$ be the ambiguity graph, and call its edges \emph{ambiguous} pairs. 
            Since for a good pair $i<j$, all vertices of $Z_j^{(2)}$ have the same type towards $i$, ambiguous pairs are bad and they are edges of $J$ (if $ij\notin E(J)$ the blocks are complete, so all vertices are complete towards $i$). 
            In particular an increasing path of $\mathcal A$ is an increasing path of $J$.

            \begin{claim}\label{clm:monotone5}
                Let $i_1<\cdots<i_p$ be an increasing path of $\mathcal A$ and $i=i_a$, $j=i_{a+1}$ be consecutive on it. 
                Let $\ell$ be a later vertex of the path. 
                Then no vertex of $Z_\ell^{(3)}$ is complete towards $j$ and sparse towards $i$.
            \end{claim}
            \begin{subproof}[Proof of Claim \ref{clm:monotone5}]
                Suppose $w\in Z_\ell^{(3)}$ is complete towards $j$ and sparse towards $i$. 
                By Definition \ref{def:frontier} (F4), some root $u\in R_j$ is anticomplete to $Z_\ell$. 
                By Definition \ref{def:frontier} (F1), $u$ is complete to $Z_j^{(3)}$.
                Since $ij$ is bad, $u$ is complete to $Z_i^{(3)}$. \
                As $Z_j^{(3)}$ contains both a vertex complete towards $i$ and a vertex sparse towards $i$ and $\G[Z_j^{(3)}]$ is connected, there is a $\G$-edge $vv'$ inside $Z_j^{(3)}$ with $v$ sparse and $v'$ complete towards $i$. 
                Since $w$ and $v$ each have fewer than $\beta|Z_i^{(3)}|$ neighbors in $Z_i^{(3)}$ and $2\beta<1$, there is $p\in Z_i^{(3)}$ nonadjacent to both $w$ and $v$. 
                Then $u-w-p-v-v'$ is an induced path of $\G$ as $uw$, $wp$, $pv$, $vv'$ are $\G$-edges, while in $G$ the vertex $u$ is complete to $\{p,v,v'\}$, the vertex $w$ is complete to $\{v,v'\}\subseteq Z_j^{(3)}$, and $v'$ is adjacent to $p$. 
                This is a rooted $P_5$ from $u$ into $Y$, contrary to Definition \ref{def:frontier} (F1).
            \end{subproof}

            By Claim \ref{clm:monotone5}, the hypothesis of Lemma~\ref{lem:monopath} holds, which gives, with $k'=\lfloor\sqrt m\rfloor$, a $k'\beta$-semisparse, hence $\rho$-semisparse, blockade of length $k'$ and width at least $\min_i|Z_i^{(3)}|/k'$. 

            Now we are ready to check the lengths and widths of all early outcome blockades. 
            Recall that $|Z_i|\ge n/m^b$ for every $i$.
            The inputs to the first and second applications of Lemma~\ref{lem:p4purify} have at least $n/m^b$ and $n/m^{b+1}$ vertices, respectively.
            Moreover, $|Z_i^{(2)}|\ge n/(2m^{b+2})$.
            Table~\ref{tab:frontier5} lists every blockade the proof can produce with its length $K$, its width as given by the lemma used and an exponent $e$ with width at least $n/K^e$. 
            Every length satisfies $m^{1/4}\le K\le1/\xi_0=64m^5/\rho\le(m/\rho)^6$, all exponents are at most $4b+8\le4b+12\le20+4b$ and $K\le(m/\rho)^{20}$. 
            This proves Lemma \ref{lem:frontier5}.
        \end{proof}

        \begin{table}[ht]
            \centering\small
            \renewcommand{\arraystretch}{1.25}
            \setlength{\tabcolsep}{4pt}
            \providecommand{\rr}{\raggedright\let\\\tabularnewline}
            \begin{tabular}{|c|l|l|l|c|p{4.1cm}|}
            \hline
            Round & Outcome & Length $K$ & Width (as given) & $e$ & \rr Why width $\ge n/K^e$ \\
            \hline
            1 & (M1) & $[m^{1/4},1/\xi_0]$ & $n/(m^bK^4)$ & $4b+4$ & \rr $m\le K^4$, so $m^b\le K^{4b}$. \\
            1 & (M2) & $[m^{3/4}/9,\,m^{3/4}/8]$ & $n/(2m^{b+1})$ & $2b+3$ & \rr $K\ge m^{3/4}/9\ge81$ gives $m\le(9K)^{4/3}\le K^2$. \\
            1 & (M3) & $\lfloor1/\xi_0\rfloor\ge m$ & $\xi_0^2n/m^{b+1}$ & $b+5$ & \rr $\xi_0\ge1/(2K)$ and $m\le K$, so width $\ge n/(4K^{b+3})$. \\
            1 & (M4) & $\lfloor1/(\xi_0m)\rfloor\ge m$ & $\xi_0^3n/(4m^b)$ & $b+11$ & \rr $\xi_0\ge1/(2mK)\ge1/(2K^2)$, so width $\ge n/(32K^{b+6})$. \\
            \hline
            2 & (M1) & $[m^{1/4},m^2]$ & $n/(m^{b+1}K^4)$ & $4b+8$ & \rr $m\le K^4$. \\
            2 & (M2) & $[m^{3/4}/9,\,m^{3/4}/8]$ & $n/(2m^{b+2})$ & $2b+5$ & \rr $m\le K^2$, as in round~1. \\
            2 & (M3) & $m^2$ & $n/m^{b+6}$ & $(b+6)/2$ & \rr Width $\xi^2n/m^{b+2}$ with $\xi=m^{-2}$ and $m=K^{1/2}$. \\
            2 & (M4) & $m$ & $n/(4m^{b+7})$ & $b+9$ & \rr Width $\xi^3n/(4m^{b+1})$ with $\xi=m^{-2}$ and $m=K$. \\
            \hline
            3 & Lemma~\ref{lem:cores} & $m$ & $n/(2m^{b+4})$ & $b+5$ & \rr Width $|Z_i^{(2)}|/m^2$; $m=K$. \\
            final & Lemma~\ref{lem:monopath} & $k'=\lfloor\sqrt m\rfloor$ & $n/(2m^{b+3}k')$ & $3b+11$ & \rr $m\le k'^3$ and $k'\ge2^8$.
            The blockade is $\rho$-semisparse since $k'\beta\le\rho$. \\
            \hline
            \end{tabular}
            \caption{Blockades produced in the proof of Lemma~\ref{lem:frontier5}. Rounds~1 and~2 apply Lemma~\ref{lem:p4purify} with parameters $(\xi_0,m)$ and $(m^{-2},m)$ to sets of at least $n/m^b$ and $n/m^{b+1}$ vertices. All blockades except the final one are pure or complete, hence $\rho$-semisparse.}
            \label{tab:frontier5}
        \end{table}

    \subsection{Frontiers of order at least six}\label{subsec:frontier-high}
        
        \begin{lemma}\label{lem:frontier-high}
            Let $r\ge6$ and $\cC$ be a hereditary class satisfying $\Rs_{r-1}$ and $\Rs_{\max(r-2,5)}$, both with constant $d_0\ge2$. 
            Then $\cC$ has the frontier property of order $r$ with constant $c=d_0(4d_0+11)$.
        \end{lemma}

        \begin{proof}
            Let $(Z_1,\ldots,Z_m)$, $R_j$, $J$ be a rooted frontier of order $r$ in $Y$, with $G\in\cC$ and $m,\rho,b,n$ as in Definition~\ref{def:FP}. 
            All applications of $\Rs_{r-1}$ and $\Rs_{\max(r-2,5)}$ below are made in $G$, with coarse parameter $1/m\le2^{-d_0}$. 
            Whenever such an application gives Definition \ref{def:Rr} (R2) or (R3), we stop (an \emph{early outcome}). 
            Put
            \[
                a_0=\frac{\rho}{2^{12}m^{4d_0+10}},\quad\theta=\frac{a_0}4,\quad a_1=\frac{\rho}{64m^2},\quad\eta=\frac1{16m^{d_0+2}},\quad\beta=2m^{d_0}\theta,\quad\beta'=2m^{d_0}\beta .
            \]

            We have $0<a_0\le a_1\le\min\{\rho,1/m\}$.
            We proceed as follows.
            
            \emph{Round 1.} By Definition \ref{def:frontier} (F2), the vertices in $\bigcup_{j>i}Z_j$ send no rooted $P_{r-1}$ into $Z_i$. 
            Since $|Z_i|\ge(m/\rho)^c\ge a_0^{-d_0}$, by applying $\Rs_{r-1}$ to $Z_i$ with $x$ replaced by $a_0$, outside early outcomes, Definition \ref{def:Rr} (R1) gives $Z_i^{(1)}\subseteq Z_i$ with $|Z_i^{(1)}|\ge|Z_i|/m^{d_0}$ such that every vertex in $\bigcup_{j>i}Z_j$ is complete to $Z_i^{(1)}$ or has fewer than $\theta|Z_i^{(1)}|$ neighbors in $Z_i^{(1)}$. 
            So, as in round~1 of Lemma~\ref{lem:frontier5}, every vertex $p\in Z_j^{(1)}$ has a \emph{type} towards each $i<j$: \emph{complete} if $p$ is complete to $Z_i^{(1)}$ and \emph{sparse} if $p$ has fewer than $\theta|Z_i^{(1)}|$ neighbors in $Z_i^{(1)}$.
            A vertex complete towards $i$ remains complete to every subset of $Z_i^{(1)}$.

            \emph{Round 2.} Call a pair $i<j$ \emph{good} if some root $u\in R_j$ has at least $\eta|Z_i^{(1)}|$ nonneighbors in $Z_i^{(1)}$ and \emph{bad} otherwise. 
            For a good pair $i<j$, fix such a root $u$ and let $S_{ij}$ be the set of vertices of $Z_j^{(1)}$ sparse towards $i$. 
            There are at most $\theta|Z_i^{(1)}||S_{ij}|$ edges between $Z_i^{(1)}$ and $S_{ij}$. 
            Since $u$ has at least $\eta|Z_i^{(1)}|$ nonneighbors in $Z_i^{(1)}$, we can fix a nonneighbor $p_{ij}$ of $u$ in $Z_i^{(1)}$ with at most $\theta|S_{ij}|/\eta\le\theta|Z_j^{(1)}|/\eta$ neighbors in $S_{ij}$. 
            We call $p_{ij}$ the \emph{representative} of the good pair $(i,j)$.
            By Definition \ref{def:frontier} (F1) and Lemma~\ref{lem:prefix}, $p_{ij}$ sends no rooted $P_{r-1}$ into $Z_j^{(1)}$. 
            Since $c\ge4d_0$, $m\ge64$ and $0<\rho\le1$, we have $|Z_j^{(1)}|\ge |Z_j|/m^{d_0}\ge(m/\rho)^c/m^{d_0}\ge(64m^2/\rho)^{d_0}=a_1^{-d_0}$.
            We may therefore apply $\Rs_{r-1}$ to $Z_j^{(1)}$ with parameters $x=a_1$ and $y=1/m$, using the vertices $p_{ij}$ as roots.
            Outside early outcomes, we obtain $D_j\subseteq Z_j^{(1)}$ with $|D_j|\ge|Z_j^{(1)}|/m^{d_0}$ such that, for every good pair $i<j$, the representative $p_{ij}$ is complete to $D_j$ or has fewer than $a_1|D_j|/4$ neighbors in $D_j$.
            As in Lemma~\ref{lem:frontier5}, the sets $D_j$ are chosen first and then we delete the following vertices. For each good pair $i<j$, delete the set $S_{ij}\cap D_j$ if $p_{ij}$ is complete to $D_j$ and delete the vertices of $D_j$ complete towards $i$ otherwise (they are neighbors of $p_{ij}$). 
            For each bad pair $j<\ell$, delete every vertex of $D_j$ that is nonadjacent to at least one root in $R_\ell$.
            For fixed $j$, the deletions of the first kind remove at most a fraction $m(a_1/4+\theta m^{d_0}/\eta)=\rho/(256m)+\rho/(2^{10}m^{2d_0+7})<1/16$ of $D_j$, and those of the second kind remove at most $m^2\eta|Z_j^{(1)}|\le|D_j|/16$ vertices.
            Let $Z_j^{(2)}$ be the remaining set, so $|Z_j^{(2)}|\ge|Z_j|/(2m^{2d_0})$.
            A vertex sparse towards $i$ has fewer than $\beta|Z_i^{(2)}|$ neighbors in $Z_i^{(2)}$. 
            For a good pair $i<j$, all vertices of $Z_j^{(2)}$ have the same type towards $i$ and for a bad pair $i<j$ every root in $R_j$ is complete to $Z_i^{(2)}$. 
            Call $i<j$ \emph{ambiguous} if $Z_j^{(2)}$ contains both a vertex complete towards $i$ and a vertex sparse towards $i$. 
            So ambiguous pairs are bad and are edges of $J$.

            \emph{Round 3.} 
            Call an ambiguous pair $i<j$ \emph{offending}, if there are $\ell>j$, $w\in Z_\ell^{(2)}$ and $u\in R_j$ such that $u$ is nonadjacent to $w$ and $w$ is complete towards $j$ and sparse towards $i$. 
            For each offending pair fix such $(\ell,w,u)$. 
            Since the pair is bad, $u$ is complete to $Z_i^{(2)}$. 
            By Definition \ref{def:frontier} (F1), $u$ is also complete to $Z_j^{(2)}$. 
            Since the vertex $w$ has fewer than $\beta|Z_i^{(2)}|$ neighbors in $Z_i^{(2)}$ and the vertices $S'_{ij}\subseteq Z_j^{(2)}$ that are sparse towards $i$ send fewer than $\beta|Z_i^{(2)}||S'_{ij}|$ edges to $Z_i^{(2)}$, some $p\in Z_i^{(2)}\setminus N(w)$ has at most $2\beta|Z_j^{(2)}|$ neighbors in $S'_{ij}$. 
            Fix such a representative $p'_{ij}=p$.

            If $p'_{ij}$ sends a rooted $P_{r-2}=p'_{ij}-q_2-\cdots -q_{r-2}$ into $Z_j^{(2)}$, then $u-w-p'_{ij}-q_2-\cdots-q_{r-2}$ would be a rooted $P_r$ from $u$ into $Y$, a contradiction. 
            So we may assume that the representative $p'_{ij}$ sends no rooted $P_{r-2}$ into $Z_j^{(2)}$.
            In particular, $p'_{ij}$ sends no rooted $P_{\max(r-2,5)}$ into $Z_j^{(2)}$. 
            Since $c\ge6d_0$, $m\ge64$, and $0<\rho\le1$, we have $|Z_j^{(2)}|\ge |Z_j|/(2m^{2d_0})\ge(m/\rho)^c/(2m^{2d_0})\ge(64m^2/\rho)^{d_0}=a_1^{-d_0}$.
            We may therefore apply $\Rs_{\max(r-2,5)}$ to $Z_j^{(2)}$ with parameters $x=a_1$ and $y=1/m$, using these (fewer than $m$) representatives as roots.
            Outside early outcomes, we get $D'_j\subseteq Z_j^{(2)}$ with $|D'_j|\ge|Z_j^{(2)}|/m^{d_0}$ on which each representative is complete or has fewer than $a_1|D'_j|/4$ neighbors. 
            For each offending pair $i<j$, we delete vertices of $D'_j$ as follows. 
            If $p'_{ij}$ is complete to $D'_j$, then delete the vertices of $D'_j$  that are sparse towards $i$, thus, at most $2\beta m^{d_0}|D'_j|$ vertices. 
            If $p'_{ij}$ has fewer than $a_1|D'_j|/4$ neighbors in $D'_j$, then delete the vertices complete towards $i$, thus, fewer than $a_1|D'_j|/4$, since they are adjacent to $p'_{ij}$. 
            For fixed $j$, there are fewer than $m$ offending pairs $i<j$, so the total number of vertices deleted from $D'_j$ is at most $m(a_1/4+2m^{d_0}\beta)|D'_j|$.
            Since $m(a_1/4+2m^{d_0}\beta)=\rho/(256m)+\rho/(2^{12}m^{2d_0+9})\le1/256+1/4096<1/2$, fewer than half the vertices of $D'_j$ are deleted.
            The remaining set $Z_j^{(3)}$ therefore satisfies $|Z_j^{(3)}|\ge|D'_j|/2\ge|Z_j^{(2)}|/(2m^{d_0})\ge|Z_j|/(4m^{3d_0})$.
            A vertex sparse towards $i$ has fewer than $\beta'|Z_i^{(3)}|<|Z_i^{(3)}|$ neighbors in $Z_i^{(3)}$.
            In particular, a vertex complete to $Z_i^{(3)}$ has type complete towards $i$.

            So the types satisfy the setting of Lemma~\ref{lem:monopath} for $W_i=Z_i^{(3)}$ and $\beta'$. 
            Let $\mathcal A'$ be the ambiguity graph. 
            Every edge of $\mathcal A'$ is ambiguous and hence an edge of $J$ and not offending.

            \begin{claim}\label{clm:monotone-high}
                Let $i_1<\cdots<i_p$ be an increasing path of $\mathcal A'$ and $i=i_a$, $j=i_{a+1}$ be consecutive on it.
                Let $\ell$ be a later vertex of the path. 
                Then no vertex of $Z_\ell^{(3)}$ is complete towards $j$ and sparse towards $i$.
            \end{claim}
            \begin{subproof}[Proof of Claim \ref{clm:monotone-high}]
                The path is also an increasing path of $J$, so by Definition \ref{def:frontier} (F4), some $u\in R_j$ is anticomplete to $Z_\ell\supseteq Z_\ell^{(3)}$. 
                If some $w\in Z_\ell^{(3)}$ were complete towards $j$ and sparse towards $i$, then $(\ell,w,u)$ would witness that $ij$ is offending, since $\ell>j$, $w\in Z_\ell^{(3)}\subseteq Z_\ell^{(2)}$, and $u\in R_j$ is nonadjacent to $w$.
                This is impossible since $ij\in E(\mathcal A')$ and no edge of $\mathcal A'$ is offending.
            \end{subproof}

            By Claim~\ref{clm:monotone-high} and Lemma~\ref{lem:monopath} with $k'=\lfloor\sqrt m\rfloor$, we obtain a $k'\beta'$-semisparse blockade of length $k'$ and width at least $\min_i|Z_i^{(3)}|/k'\ge n/(4m^{b+3d_0}k')$.
            Since $k'\le m$ and $\beta'=\rho/(2^{12}m^{2d_0+10})$, we have $k'\beta'\le\rho$, so this blockade is $\rho$-semisparse.
            
            It remains to verify that each possible outcome satisfies Definition~\ref{def:FP}.
            The bounds above imply that the input set of each application has at least $n/(2m^{b+2d_0})$ vertices.
            An early (R2) outcome gives a $z^4$-restricted set $T$ with $z\in[a_0^{d_0},1/m]\subseteq[(\rho/m)^c,1/m]$.
            Its size satisfies $|T|\ge z^{d_0}n/(2m^{b+2d_0})\ge z^{b+3d_0+1}n\ge z^{c+b}n$.
            An early (R3) outcome gives a $\rho$-semisparse blockade whose length satisfies $K\in[m,a_0^{-d_0}]\subseteq[m^{1/4},(m/\rho)^c]$.
            Its width is at least $n/(2m^{b+2d_0}K^{d_0})\ge n/K^{b+3d_0+1}\ge n/K^{c+4b}$.
            The final blockade has length $K=k'=\lfloor\sqrt m\rfloor$, with $m^{1/4}\le K\le(m/\rho)^c$.
            Since $m\le K^3$ and $K\ge4$, its width is at least $n/(4m^{b+3d_0}K)\ge n/(4K^{3b+9d_0+1})\ge n/K^{3b+9d_0+2}\ge n/K^{c+4b}$.
            This proves Lemma \ref{lem:frontier-high}.
        \end{proof}

    \subsection{Proof of Theorems~\ref{thm:R5} and~\ref{thm:Rstep}}\label{subsec:Rproof}

        \begin{theorem}\label{thm:purify}
            Let $r\ge5$ and let $\cC$ be a hereditary class with the frontier property of order $r$ with constant $c$, such that $\cO_{r-1}(\cC)$ has Erd\H{o}s--Hajnal exponent $\kappa\in(0,1]$. 
            Let $g=\max(2,1/\kappa)$. 
            Then $\cC$ satisfies $\Rs_r$ with constant $d=(c+6)\bigl(6\lceil4g+5\rceil+3\bigr)+c+1$.
        \end{theorem}
        
        \begin{proof}
            The proof combines the cut-tree lemma, the frontier property, and the sampling lemma.
            We choose the parameters so that each outcome of the cut-tree lemma yields one of Definition \ref{def:Rr} (R1)--(R3) with the bounds required by $\Rs_r$.
            
            Put $d_t=\lceil4g+5\rceil$ and $d_1=6d_t+3$. 
            So $d=(c+6)d_1+c+1$. 
            Let $x,y,G,U,Y,n$ be as in Definition~\ref{def:Rr}, and put $t=\lceil y^{-d_t}\rceil$. 
            Then $t\ge2^{d\,d_t}\ge\max(2^{16},2^c)$ and $t\le2y^{-d_t}\le x^{-d_1}$. 
            
            We then apply Lemma~\ref{lem:tree} with $\mu=n/t^6$ and show that each of its outcomes yields one of (R1)--(R3) in Definition~\ref{def:Rr}.

            (T1) gives (R1): the set $S$ is $x$-clean for $U$, and since $t\le2y^{-d_t}$, $2^{-6}\ge y$ and $d\ge6d_t+1$, we have $|S|\ge n/t^6\ge2^{-6}y^{6d_t}n\ge y^{6d_t+1}n\ge y^dn$.

            (T2) gives (R3) with $K=t\in[1/y,x^{-d}]$: a complete blockade is $x$-semisparse, and $n/t^7\ge n/K^d$.

            (T3) and (T4) give a rooted frontier with $m$ blocks, where $b=6$ or $b=3$ respectively, and $t\le m\le\max(t,4t^6/x)$. 
            Since $y\ge x$ and $2^8\le x^{-2}$, $4t^6/x\le2^8y^{-6d_t}x^{-1}\le x^{-6d_t-3}$, and so $m\leq x^{-6d_t-3}= x^{-d_1}$. 
            We have $m\ge t\ge2^c$, and every block has at least $n/m^b$ vertices.
            Since $b\le6$, $n\ge x^{-d}$, $m\le x^{-d_1}$, and $d-6d_1=c(d_1+1)+1$, we have
            \[
                \frac{n}{m^b}\ge\frac{n}{m^6}\ge x^{-d+6d_1}=x^{-c(d_1+1)-1}\ge x^{-c(d_1+1)}\ge(m/x)^c.
            \]
            Thus the frontier property applies with $\rho=x$.
            If the first alternative holds, there are $z\in[(x/m)^c,1/m]$ and a $z^4$-restricted set $T\subseteq Y$ with $|T|\ge z^{c+b}n$.
            Since $d\ge c(d_1+1)$, the bounds on $m$ give $z\ge(x/m)^c\ge x^{c(d_1+1)}\ge x^d$ and $z\le1/m\le1/t\le y$.
            Moreover, $b\le6$ and $d\ge c+6$ imply $|T|\ge z^{c+b}n\ge z^{c+6}n\ge z^dn$, since $z\le1$.
            Hence (R2) holds.
            If the second alternative holds, there is an $x$-semisparse blockade of length $K\in[m^{1/4},(m/x)^c]$ and width at least $n/K^{c+4b}$.
            Since $d_t\ge4$ and $d\ge c(d_1+1)$, we have
            \[
                \frac1y\le y^{-d_t/4}\le t^{1/4}\le m^{1/4}\le K\le(m/x)^c\le x^{-c(d_1+1)}\le x^{-d}.
            \]
            Its width is at least $n/K^{c+4b}\ge n/K^{c+24}\ge n/K^d$, since $b\le6$ and $d\ge c+24$.
            Hence (R3) holds.

            (T5): let $M\ge3n/(8t)$ be the total size of the negative layers. 
            Each layer has fewer than $n/t^6$ vertices, so $\varphi\le8/(3t^5)\le\frac83y^{5d_t}\le y^{16g+16}$. By (T5), every set $X$ meeting each layer at most once induces a member of $\cO_{r-1}(\cC)$, so $\homog(G[X])\ge|X|^\kappa$. 
            By Lemma~\ref{lem:sampling} with $\eps=y^4\le2^{-4(g+10)}$, there is a $y^4$-restricted set of size at least $y^{8g+8}M\ge y^{8g+d_t+10}n\ge y^dn$, which is (R2) with $z=y$.

            Finally $d\ge\max(6d_t+2,\,c(d_1+1),\,c+24,\,8g+d_t+10,\,g+10)$, which covers all the exponents used.
        \end{proof}

        Now we are ready to prove Theorems \ref{thm:R5} and \ref{thm:Rstep}, which we restate here for convenience. 
        
        \begin{theorem}\label{thm:R5 2}\textnormal{(Theorem \ref{thm:R5})}
            Every hereditary class satisfies $\Rs_5$ with constant $2127$.
        \end{theorem}
        
        \begin{proof}
            Let $\cC$ be any hereditary class. We claim that $\cO_4(\cC)$ has Erd\H{o}s--Hajnal exponent $\frac12$.
            Let $G\in\cO_4(\cC)$ with ordering $v_1,\ldots,v_n$. 
            Suppose that $a-b-c-d$ is an induced path of $\overline{G}$.
            If the earliest vertex of this $P_4$ is $a$ or $d$, then this vertex would start a rooted $P_4$ into the later vertices. 
            So the earliest vertex is $b$ or $c$. 
            Hence $b<a$ or $c<d$, which satisfies the hypothesis of Theorem~\ref{thm:chvatal} for $\overline{G}$.
            Thus $\homog(G)=\homog(\overline{G})\ge|G|^{1/2}$.
            By Lemma~\ref{lem:frontier5} and Theorem~\ref{thm:purify} with $c=20$ and $\kappa=\frac12$ (so $g=2$, $d_t=13$, and $d_1=81$), $\cC$ satisfies $\Rs_5$ with constant $26\cdot81+21=2127$.
        \end{proof}

        \begin{theorem}\label{thm:Rstep 2}\textnormal{(Theorem \ref{thm:Rstep})}
            Let $r\ge6$ and let $\cC$ be a hereditary class. Suppose that $\cC$ satisfies $\Rs_{r-1}$ and $\Rs_{\max(r-2,5)}$, and that $\cO_{r-1}(\cC)$ has the Erd\H{o}s--Hajnal property. Then $\cC$ satisfies $\Rs_r$.
        \end{theorem}

        \begin{proof}
            Let $d_0\ge2$ be a common constant for $\Rs_{r-1}$ and $\Rs_{\max(r-2,5)}$. By Lemma~\ref{lem:frontier-high}, $\cC$ has the frontier property of order $r$ with constant $d_0(4d_0+11)$. 
            By Theorem~\ref{thm:purify}, $\cC$ satisfies $\Rs_r$.
        \end{proof}

\section{Proof of the main theorem}\label{sec:proof}

    Throughout this section, fix $s\ge 6$ and assume $\EH(P_{s-1})$. 
    Let $\cG_s$ be the class of $\overline{P_s}$-free graphs.

    \begin{lemma}\label{lem:ordered}
        Let $r\ge5$. If $\cG_s$ satisfies $\Rs_r$, then $\cO_r(\cG_s)$ has the local alternative.
    \end{lemma}
    \begin{proof}
        Let $d$ be the constant of $\Rs_r$. 
        We show that $\cO_r(\cG_s)$ has the local alternative with constant $D=\max(d+2,4)$. 
        Let $0<x\le y\le2^{-D}$ and $F\in\cO_r(\cG_s)$ be $y^3$-sparse with $n=|F|\ge x^{-D}$. 
        Order $V(F)$ as in Definition~\ref{def:Er}, let $Y$ be the set of the last $\lceil yn\rceil$ vertices and $U$ the set of the others. 
        Every vertex of $U$ precedes all of $Y$, so it sends no rooted $P_r$ into $Y$. 
        Since $n\ge x^{-D}\ge1/y$ we have $yn\ge1$.
        So $yn\le|Y|\le yn+1\le2yn$. 
        Since $D-1\ge d$, $|Y|\ge yn\ge xn\ge x^{1-D}\ge x^{-d}$. 
        Also $x\le y\le2^{-D}\le2^{-d}$. 
        So the conditions of Definition \ref{def:Rr} hold for $F$ with $x,y$. 

        If Definition \ref{def:Rr} (R1) holds, then there is a set $X\subseteq Y$ that is $x$-clean for $U$ and has $|X|\ge y^d|Y|\ge y^{d+1}n$. 
        Fix a vertex $p\in X$. 
        Every vertex complete to $X$ is a neighbor of $p$, and $p$ has at most $y^3n$ neighbors as $F$ is $y^3$-sparse. 
        So at most $y^3n$ vertices of $U$ are complete to $X$. 
        Since $X$ is $x$-clean for $U$, the other vertices of $U$ have fewer than $x|X|/4$ neighbors in $X$. 
        Let $Z=\{u\in U:|N_F(u)\cap X|<x|X|/4\}$.
        Then $X$ and $Z$ are disjoint, $Z$ is $x$-sparse to $X$, and $|Z|\ge|U|-y^3n\ge(1-4y)n$.
        Since $|X|\ge y^{d+1}n\ge y^Dn$, Definition \ref{def:LA} (O3) holds.

        If Definition \ref{def:Rr} (R2) holds, then there is a $z^4$-restricted set with $z\in[x^d,y]$ and at least $z^d|Y|\ge z^dyn\ge z^{d+1}n$ vertices (as $z\le y$), which gives Definition \ref{def:LA} (O4).

        If Definition \ref{def:Rr} (R3) holds, then there is an $x$-semisparse (hence weakly $x$-semisparse) blockade of length $K\in[1/y,x^{-d}]$ and width at least $|Y|/K^d\ge yn/K^d\ge n/K^{d+1}$ (as $y\ge1/K$), which gives Definition \ref{def:LA} (O2).
            
        Thus $\cO_r(\cG_s)$ has the local alternative. 
    \end{proof}

    The next lemma takes the place of property $(*)$ of \cite[Section~4]{hjz} and of the quantitative comb condition \cite[Definition~7.1]{tn}.
        
    \begin{lemma}\label{lem:combstep}
        If $\cG_s$ satisfies $\Rs_{\max(s-2,5)}$, then $\cG_s$ has the local alternative.
    \end{lemma}
    \begin{proof}
        Let $r=\max(s-2,5)$ and $d$ be the constant of $\Rs_r$.
        We show $\cG_s$ has the local alternative with constant $D=3d+10$. 
        Let $0<x\le y\le2^{-D}$ and $G\in\cG_s$ be $y^3$-sparse with $n=|G|\ge x^{-D}$. 
        By applying Lemma~\ref{lem:comb} with these $x,y$, there is an $(\ell,y^4n/\ell^2)$-comb $((a_i,B_i):i\in[\ell])$ with $\ell\in[1/y,x^{-2}]$ and a special vertex $v$, since the first two outcomes of Lemma~\ref{lem:comb} give Definition \ref{def:LA} (O3) and (O1). 
        As $y\ge1/\ell$, each $|B_i|\ge n/\ell^6$.

        Let $U_i=\bigcup_{j\ne i}B_j$.
        In $\G$, the vertex $v$ is adjacent to $a_i$ and to no vertex of $B_i\cup U_i$, while $a_i$ is adjacent to every vertex of $U_i$ and to no vertex of $B_i$.
        If some $u\in U_i$ sent a rooted $P_{s-2}$ into $B_i$, then prefixing this path with $v,a_i$ would give an induced $P_s$ in $\G$, contrary to the fact that $\G$ is $P_s$-free.
        Therefore no vertex of $U_i$ sends a rooted $P_{s-2}$ into $B_i$. 
        Since $r\ge s-2$, no vertex of $U_i$ sends a rooted $P_r$ into $B_i$.

        Since $x^3\le1/\ell\le y\le2^{-d}$ and $|B_i|\ge y^4n/\ell^2\ge x^8n\ge x^{-3d}$, the conditions of $\Rs_{r}$ hold for $Y=B_i$ and the roots $U=U_i$ with $x,y$ replaced by $x^3,1/\ell$. 
        If Definition \ref{def:Rr} (R2) holds, then there is a $z^4$-restricted set with $z\in[x^{3d},1/\ell]$ and at least $z^d|B_i|\ge z^{d+6}n$ vertices, which gives Definition \ref{def:LA} (O4).
        If Definition \ref{def:Rr} (R3) holds, then there is an $x^3$-semisparse blockade of length $K\in[\ell,x^{-3d}]$ and width at least $|B_i|/K^d\ge n/K^{d+6}$, and hence Definition \ref{def:LA} (O2) holds, since $K\ge\ell\ge1/y$. 
    
        So we may assume that Definition \ref{def:Rr} (R1) holds for every $i$, giving $S_i\subseteq B_i$ with $|S_i|\ge|B_i|/\ell^d\ge n/\ell^{d+6}$ that is $x^3$-clean for $U_i$. 
        We claim that $(S_1,\ldots,S_\ell)$ is a weakly $x$-semisparse blockade. 
        Let $\theta=x^3/4$ and $i\ne j$. 
        Since $S_j\subseteq B_j\subseteq U_i$, every vertex of $S_j$ is a root for $S_i$. 
        So by the cleanliness of $S_i$, every vertex of $S_j$ is complete to $S_i$ or has fewer than $\theta|S_i|$ neighbors in $S_i$.
        Symmetrically, since $S_i\subseteq U_j$, every vertex of $S_i$ is complete to $S_j$ or has fewer than $\theta|S_j|$ neighbors in $S_j$. 
        Suppose $(S_i,S_j)$ is not complete. 
        Take nonadjacent $p\in S_i$, $q\in S_j$. 
        Then $q$ is not complete to $S_i$, so it has fewer than $\theta|S_i|$ neighbors in $S_i$. 
        In particular, fewer than $\theta|S_i|$ vertices of $S_i$ are complete to $S_j$. 
        Every other vertex of $S_i$ is not complete to $S_j$, so it has fewer than $\theta|S_j|$ neighbors in $S_j$. 
        Hence there are fewer than $2\theta|S_i||S_j|\le x|S_i||S_j|$ edges between $S_i$ and $S_j$. 
        This gives Definition \ref{def:LA} (O2) with $k=\ell$.

        All exponents above are at most $3d+10$, which proves the lemma.
    \end{proof}

    Now we are ready to prove Theorem \ref{thm:main}. 

    \begin{proof}[Proof of Theorem~\ref{thm:main}]
        We prove by induction on $s$.
        For $s\le4$ the result is classical and $\EH(P_5)$ follows from \cite{nss7}. 
        So we may assume that $\EH(P_{s-1})$. 

        We claim that $\cG_s$ satisfies $\Rs_r$ for $5\le r\le \max\{s-2,5\}$.
        We prove this by induction on $r$. 
        For $r=5$, this is Theorem~\ref{thm:R5}. 
        Let $6\le r\le s-2$ and suppose $\cG_s$ satisfies $\Rs_{r-1}$ and $\Rs_{\max(r-2,5)}$.
        By Lemma~\ref{lem:ordered}, $\cO_{r-1}(\cG_s)$ has the local alternative. 
        Since $\cO_{r-1}(\cG_s)$ is a hereditary class of $\overline{P_s}$-free graphs, $\cO_{r-1}(\cG_s)$ has the Erd\H{o}s--Hajnal property by Theorem~\ref{thm:machinery}.
        So Theorem~\ref{thm:Rstep} shows that $\cG_s$ satisfies $\Rs_r$. 
        Thus, $\cG_s$ satisfies $\Rs_{\max\{s-2,5\}}$. 
        By Lemma~\ref{lem:combstep} and Theorem~\ref{thm:machinery}, $\cG_s$ has the Erd\H{o}s--Hajnal property.
        Hence $\EH(\overline{P_s})$ holds, which is equivalent to $\EH(P_s)$.
    \end{proof}

\section*{Acknowledgements}

    The proof presented in this paper was developed through an interactive dialogue with Claude Opus 5.5 (Anthropic) and ChatGPT-6 (OpenAI). The authors used them to explore multiple possible proof strategies and directions, and the AI generated the initial draft under the authors' guidance. The authors then critically examined the generated proof, verified each step for correctness and rigor, and revised it where necessary. The authors take full responsibility for the final content.


\begin{thebibliography}{99}\small
\bibitem{aps} N.~Alon, J.~Pach and J.~Solymosi. Ramsey-type theorems with forbidden subgraphs. \emph{Combinatorica} 21 (2001), 155--170.

\bibitem{chvatal} V.~Chv\'atal. Perfectly ordered graphs. In \emph{Topics on Perfect Graphs}, Annals of Discrete Mathematics 21 (1984), 63--65.

\bibitem{chudnovsky} M.~Chudnovsky. The Erd\H{o}s--Hajnal conjecture---a survey. \emph{J. Graph Theory} 75 (2014), 178--190.

\bibitem{cs} M.~Chudnovsky and S.~Safra. The Erd\H{o}s--Hajnal conjecture for bull-free graphs. \emph{J. Combin. Theory Ser. B} 98 (2008), 1301--1310.

\bibitem{csss} M.~Chudnovsky, A.~Scott, P.~Seymour and S.~Spirkl. Erd\H{o}s--Hajnal for graphs with no 5-hole. \emph{Proc. Lond. Math. Soc.} 126 (2023), 997--1014.

\bibitem{Er47}
P.~Erd\H{o}s.
\newblock Some remarks on the theory of graphs.
\newblock {\em Bulletin of the American Mathematical Society}, 53:292-294, 1947.

\bibitem{eh} P.~Erd\H{o}s and A.~Hajnal. Ramsey-type theorems. \emph{Discrete Appl. Math.} 25 (1989), 37--52.

\bibitem{gallai} T.~Gallai. Transitiv orientierbare Graphen. \emph{Acta Math. Acad. Sci. Hungar.} 18 (1967), 25--66.

\bibitem{hjz} 
S.~Huang, Y.~Ju and Y.~Zhou. 
\newblock Erd\H{o}s--Hajnal conjecture beyond five-vertex graphs. 
\newblock arXiv:2606.06258v3, 2026.

\bibitem{nguyen26} T.~H.~Nguyen. Erd\H{o}s--Hajnal for excluding paths. arXiv:2610.24887, 2026.
\bibitem{nss5} T.~Nguyen, A.~Scott and P.~Seymour. Induced subgraph density. V. All paths approach Erd\H{o}s--Hajnal. arXiv:2307.15032, 2023.
\bibitem{nss7} T.~Nguyen, A.~Scott and P.~Seymour. Induced subgraph density. VII. The five-vertex path. \emph{Proc. Lond. Math. Soc.} 132(3) (2026), e70133. arXiv:2312.15333.

\bibitem{Ra30}
F.~P.~Ramsey.
\newblock On a problem of formal logic.
\newblock {\em Proceedings of the London Mathematical Society}, 30:264-286, 1930.

\bibitem{rodl} V.~R\"odl. On universality of graphs with uniformly distributed edges. \emph{Discrete Math.} 59 (1986), 125--134.
\bibitem{tn} V.-H.~Tran and T.~M.~Nguyen. The Erd\H{o}s--Hajnal property for the six-vertex graph with edge set $\{ab,bc,cd,de,af,bf,df\}$. arXiv:2608.28551, 2026.
\bibitem{zenodo} Y.~Zhou. The Erd\H{o}s--Hajnal property for paths. Zenodo, 9 October 2026, \href{https://doi.org/10.5281/zenodo.23270286}{doi:10.5281/zenodo.23270286}. Repository (public since 7 October 2026) with the Lean formalization: \url{https://github.com/OctopusTide/erdos-hajnal-all-paths}.
\end{thebibliography}
\end{document}